\documentclass{article}

\usepackage{amstext,amsmath,amssymb}
\usepackage{epsfig,rotating,multirow}
\usepackage{algorithm,algorithmic}
\usepackage{comment}
\usepackage{geometry}
\usepackage{booktabs}
\usepackage{placeins}

\usepackage{graphicx,bm,xcolor}
\usepackage{subcaption}
\usepackage{amsmath, amssymb, amsthm}
\usepackage{ulem}

\newtheorem{theorem}{Theorem}[section]
\newtheorem{lemma}[theorem]{Lemma}
\newtheorem{corollary}[theorem]{Corollary}

\numberwithin{equation}{section}

\def\CC{{\mathchoice {\setbox0=\hbox{$\displaystyle\rm C$}\hbox{\hbox
to0pt{\kern0.4\wd0\vrule height0.9\ht0\hss}\box0}}
{\setbox0=\hbox{$\textstyle\rm C$}\hbox{\hbox
to0pt{\kern0.4\wd0\vrule height0.9\ht0\hss}\box0}}
{\setbox0=\hbox{$\scriptstyle\rm C$}\hbox{\hbox
to0pt{\kern0.4\wd0\vrule height0.9\ht0\hss}\box0}}
{\setbox0=\hbox{$\scriptscriptstyle\rm C$}\hbox{\hbox
to0pt{\kern0.4\wd0\vrule height0.9\ht0\hss}\box0}}}}

\def\ZZ{{\mathchoice {\hbox{$\sf\textstyle Z\kern-0.4em Z$}}
{\hbox{$\sf\textstyle Z\kern-0.4em Z$}}
{\hbox{$\sf\scriptstyle Z\kern-0.3em Z$}}
{\hbox{$\sf\scriptscriptstyle Z\kern-0.2em Z$}}}}

\newcommand{\lw}[1]{\smash{\lower2.ex\hbox{#1}}}

\definecolor{gray}{gray}{0.6}

\newcommand{\vertiii}[1]{{\left\vert\kern-0.25ex\left\vert\kern-0.25ex\left\vert #1 
		\right\vert\kern-0.25ex\right\vert\kern-0.25ex\right\vert}}

\def\wtb{\widetilde{b}}

\def\wtt{\widetilde{t}}

\def\wtT{\widetilde{T}}

\begin{document}

\title{Explicit inverse of symmetric, tridiagonal near Toeplitz matrices Part II: with weakly diagonally dominant Toeplitz}

\author{Bakytzhan~Kurmanbek\thanks{Nazarbayev University, Department of Mathematics, 53 Kabanbay Batyr Ave, Astana 010000, Kazakhstan (bakytzhan.kurmanbek@nu.edu.kz)} \and Yogi~Erlangga\thanks{Zayed University, Department of Mathematics, Abu Dhabi Campus, P.O. Box 144534, United Arab Emirates (yogi.erlangga@zu.ac.ae)} \and Yerlan~Amanbek\thanks{Nazarbayev University, Department of Mathematics, 53 Kabanbay Batyr Ave, Astana 010000, Kazakhstan (yerlan.amanbek@nu.edu.kz)}}
\maketitle
\begin{abstract}
In this paper, we provide explicit formulas for the exact inverses of the symmetric tridiagonal near-Toeplitz matrices characterized by weak diagonal dominance in the Toeplitz part. 
Furthermore, these findings extend to scenarios where the corners of the near Toeplitz matrices lack diagonal dominance ($|\widetilde{b}| < 1$). 
Additionally, we compute the row sums and traces of the inverse matrices, thereby deriving upper bounds for their infinite norms. To demonstrate the practical applicability of our theoretical results, we present numerical examples addressing numerical solution of the Fisher problem using the fixed point method. Our findings reveal that the convergence rates of fixed-point iterations closely align with the expected rates, and there is minimal disparity between the upper bounds and the infinite norm of the inverse matrix. Specifically, this observation holds true for $|b| = 2$ with $|\widetilde{b}| \geq 1$. In other cases, there exists potential to enhance the obtained upper bounds.
\end{abstract}


\newcommand{\bs}[1]{\boldsymbol{#1}}

\section{Introduction} \hfill\\
In part 1 of this paper \cite{kurmanbek2024explicit}, we showed the importance of the banded near-Toeplitz matrices for various fields. Specifically, in solving differential equations using a finite difference method may result to the system of equations which is involved a banded near-Toeplitz matrix. 
In nonlinear problems, the fixed point iteration is used to study the convergence of numerical solution that requires a bound of inverse matrices. To provide an application of the problem, we repeat the differential equation from part 1 of this paper \cite{kurmanbek2024explicit}.

Let $\displaystyle  \frac{d^2 u}{dx^2} = f(u),  \quad  x \in \Omega = (0,L)$, and boundary conditions.
For example, $f(u)=ku(1-u)$ for the fisher equation, 
$f(u)=ke^u$ for the Gelfand-Bratu equation.

Usage of the finite difference method can lead to
 \begin{equation}
      \frac{d^2 u}{dx^2} (x_i) \approx \frac{1}{h^2} \left( u_{i-1} - 2u_{i}+u_{i+1}\right),
\end{equation}

where $x_i = ih$ and $u_i \equiv u(x_i)$. As a result, we obtain the system of nonlinear equations as follows,

\begin{equation}
    \widehat{T}_n\mathbf{u}=h^2f(\mathbf{u})
\end{equation}
where $\boldsymbol{u} = (u_1, \dots, u_n)^T \in \mathbb{R}^n$.
In the fixed-point iteration analysis, we consider the iteration
\begin{equation}
    \mathbf{u}^{k+1}=h^2 \widehat{T}^{-1}_n f(\mathbf{u}^{k}).
\end{equation}

Thus, the error of $k^{th}$ iteration is as follows,
\begin{align*}
    \| \mathbf{u}^{k+1}-\mathbf{u}^{k}\|=\|h^2 \widehat{T}^{-1}_n \left( f(\mathbf{u}^{k})-f(\mathbf{u}^{k-1})\right)\| \le h^2 \| \widehat{T}^{-1}_n  \| L_c \| \mathbf{u}^{k}-\mathbf{u}^{k-1}\|
\end{align*}

where $L_c$ is the Lipschitz constant, depends on properties of $f(\mathbf{u})$. The condition $h^2 \| \widehat{T}^{-1}_n  \| L_c < 1$ allows to have convergence in the fixed point algorithm.
Thus, the bound of $\| \widehat{T}^{-1}_n  \|$ is useful in the numerical analysis.
Several norm bounds of the inverses of some matrices have been studied for various applications \cite{dow2002explicit, heinig2013algebraic, huang1997analytical, lewis1982inversion, schlegel1970explicit, usmani1994inversion, wang2015explicit, yamamoto1979inversion,kurmanbek2021explicit,kurmanbek2021inverse}.

Consider the symmetric tridiagonal near-Toeplitz matrix
\begin{equation}
   \widehat{T}_n = \begin{bmatrix}
             \overline{b} & \hat{c} &   & &\\
             \hat{c} & \hat{b} &  \ddots & &\\
               &  \ddots   &  \ddots & \ddots & \\
               &              &  \ddots   & \hat{b} & \hat{c} \\
               &              &  &  \hat{c}  & \overline{b}
            \end{bmatrix} = -\hat{c}
            \underbrace{ \begin{bmatrix}
              \widetilde{b} & -1 &   & &\\
             -1 & b &  \ddots & &\\
               &  \ddots   &  \ddots & \ddots & \\
               &              &  \ddots  & b & -1 \\
               &              &   & -1  & \widetilde{b} 
            \end{bmatrix}}_{\widetilde{T}_n} =: -\hat{c} \, \widetilde{T}_n, \label{eq:matrixTn}
\end{equation}
where $n > 3$, $b = -\hat{b}/\hat{c}$, and $\widetilde{b} = -\overline{b}/\hat{c}$. One can express the near-Toeplitz matrix $\wtT_n$ as rank-2 decomposition
\begin{equation}
    \widetilde{T}_n = T_n + (\widetilde{b}-b) UU^T,
\end{equation}
where $T_n = \text{tridiag}(-1,b,-1)$ is the Toeplitz part of $\widetilde{T}_n$ and the matrix $U^T$ is represented by
\begin{equation}
   U^T = \begin{bmatrix}
                  1 & 0 & \dots & 0 & 0 \\
                  0 & 0 & \dots & 0 & 1
              \end{bmatrix}.
\end{equation}
The inverse of the Toeplitz part $T_n$, assuming it exits, can be defined and represented by the roots of the polynomial
$p(r) = -r^2 + br - 1$ (as shown in Dow \cite{dow2002explicit}). If $|b| \ge  2$, then $p(r)$ has two distinct real or equal roots, with the latter occurring when $b = 2$. Otherwise, $p(r)$ has two complex roots.

Suppose that $\wtT_n$ has an inverse.
By applying the Sherman-Morrison formula, the inverse of $\wtT_n$ is given by
\begin{equation}
  \widetilde{T}^{-1}_n = T^{-1}_n - (\widetilde{b} - b) T^{-1}_n U M^{-1} U^T T^{-1}_n, \label{eq:invTtilde}
\end{equation}
where $M = I_2 + (\widetilde{b}-b) U^T T^{-1}_{n} U \in \mathbb{R}^{2\times 2}$. Let $\widetilde{T}^{-1}_n = [\wtt^{-1}_{i,j}]$, $T^{-1}_n = [t^{-1}_{i,j}]$, and $\beta = \wtb - b$. 
The elements of the 2-by-2 matrix $M$ are given explicitly by
\begin{equation}
\begin{aligned}
   m_{11} = m_{22} &= 1 + (\widetilde{b} - b) t^{-1}_{1,1} = 1 + \beta t^{-1}_{n, n}, \\
   m_{12} = m_{21} &= (\widetilde{b}-b) t^{-1}_{1,n} = \beta t^{-1}_{n, 1}, 
   \end{aligned} \label{eq:elementM}
\end{equation}
obtained using the centrosymmetry and the symmetry of $T_n$. The $(i,j)$-entry of the inverse of  $\wtT_n$, denoted by $\wtt^{-1}_{i,j}$, reads
\begin{equation}
 \wtt_{i,j}^{-1} =t^{-1}_{i,j} - \frac{\beta}{\Delta} \left[ t_{i,1}^{-1}(m_{11} t^{-1}_{1,j} - m_{12} t^{-1}_{n,j}) + t^{-1}_{i,n}(-m_{12} t^{-1}_{1,j} + m_{11} t^{-1}_{n,j})\right],  \label{eq:wtt}
\end{equation}
where $\Delta = m_{11}^2 - m_{12}^2 = (1 + \beta t^{-1}_{n,n} )^2 - \beta^2  (t^{-1}_{n,1})^2$ is the determinant of $M$.

With the inverse of $M$ known explicitly, the trace of $\widetilde{T}^{-1}_n$ can be also calculated and is given by the formula
\begin{equation}
  \text{Tr}(\widetilde{T}^{-1}_n) := \sum_{i=1}^n \wtt^{-1}_{i,i} = \sum_{i=1}^n t^{-1}_{i,i} - \frac{\beta}{\Delta} \left\{ m_{11} \sum_{i=1}^n \left( (t^{-1}_{i,1})^2 + (t^{-1}_{n,i})^2 \right) - 2 m_{12} \sum_{i=1}^n t^{-1}_{i,1} t^{-1}_{n,i} \right\}.
  \label{eq:traceTt}
\end{equation}
Furthermore, the sum of the elements in the $i$-th row of the inverse can be derived:
\begin{align}
   \sum_{j=1}^n \wtt^{-1}_{i,j} &= \sum_{j=1}^n t^{-1}_{i,j} - \frac{\beta}{\Delta}\left[ t_{i,1}^{-1}\left( m_{11}\sum_{j=1}^n  t^{-1}_{1,j} - m_{12} \sum_{j=1}^n  t^{-1}_{n,j} \right) + t^{-1}_{i,n}\left( -m_{12} \sum_{j=1}^n  t^{-1}_{1,j} + m_{11} \sum_{j=1}^n  t^{-1}_{n,j}\right) \right] \notag \\
&= \sum_{j=1}^n t^{-1}_{i,j} - \frac{\beta}{m_{11} + m_{12}} (t^{-1}_{i,1} + t^{-1}_{i,n}) \sum_{j=1}^n t^{-1}_{n,j}, \label{eq:rowsumTt}
\end{align}
after making use of the centrosymmetry of $T_n$ and $T^{-1}_n$, with
$$
  m_{11} + m_{12} = 1 + \beta(t^{-1}_{n,1} + t^{-1}_{n,n}).
$$
We refer the reader to~\cite{kurmanbek2024explicit} for detailed derivations of~\eqref{eq:elementM}--\eqref{eq:rowsumTt}.

In~\cite{kurmanbek2024explicit}, we presented properties of the inverse of $\wtT$ in~\eqref{eq:matrixTn}, where $|b| > 2$, namely where the Toeplitz part is strictly diagonally dominant. In this current paper, we derive inverse properties of $\wtT$ when $|b| = 2$, i.e., when the Toeplitz part is only {\it weakly} diagonally dominant. 

In this paper, our objective is to derive the trace \eqref{eq:traceTt} and the rowsum \eqref{eq:rowsumTt}, along with establishing bounds for infinite norms of the explicit inverse of the near Toeplitz matrix $\widetilde{T}_n$ in~\eqref{eq:matrixTn}. Our results will cover both dominant ($|\widetilde{b}| \geq 1$) and non-dominant ($|\widetilde{b}| < 1$) scenarios at the corners of these near Toeplitz matrices, with a focus on $|b| = 2$. Numerical examples are provided for the Fisher problem.
\\
In Section~\ref{sec:prelim_results}, we present common preliminary results for two cases. Subsequently, in Section~\ref{sec:b=2}, we delve into the specific results when $b = 2$. In Section~\ref{sec:b=-2}, we extend the work from the previous sections to the case when $b = -2$. Additionally, Section~\ref{sec:num} provides numerical experiments addressing Fisher's problem.


\section{Preliminary results}\label{sec:prelim_results}

In this section, we summarize some preliminary results, which are useful in deriving the properties of $\wtT^{-1}_n$. We begin with the result of~\cite{dow2002explicit} regarding elements of the inverse of the symmetric tridiagonal Toeplitz matrix $T_n$ when $b = \pm 2$: for $i \geq j$,
\begin{equation} 
    t^{-1}_{i, j} = \left(\frac{2}{b}\right)^{i+1-j} \frac{j(n+1-i)}{n+1}. \label{tij_general}
\end{equation} 

\begin{lemma}\label{lemma:prem1}
Let $\wtT_{n}$ be invertible with $b = \pm 2$. Then the elements of the inverse matrix $\wtT_n^{-1}$ are, for $i \geq j$,
    \begin{equation}
        \wtt^{-1}_{i,j} = \left( \frac{2}{b} \right)^{i+1-j} \frac{\left(j \left(1 + \frac{2\beta}{b}\right) - \frac{2\beta}{b}\right) \left((n-i) \left(1 + \frac{2\beta}{b}\right) + 1 \right)}{\left(1 + \frac{2\beta}{b}\right) \left( n+1 + \frac{2\beta}{b}(n-1)\right)}, \label{gen_wtt}
    \end{equation}
    where $\beta = \wtb - b$, and $\wtt^{-1}_{i,j} = \wtt^{-1}_{j,i}$ for $i < j$.
\end{lemma}
\begin{proof}
    We first note that, using~\eqref{tij_general}, the equations in ~\eqref{eq:elementM} become
    \begin{align*}
        m_{11} = 1 + \frac{2\beta}{b} \frac{n}{n+1} \quad \text{and}  \quad m_{12} = \beta \left( \frac{2}{b}\right)^{n} \frac{1}{n+1}.
    \end{align*}
    Using the above expressions and \eqref{tij_general}, the two sums in the right-hand side of \eqref{eq:wtt} can be calculated:
    \begin{align*}
        m_{11} t^{-1}_{1,j} - m_{12} t^{-1}_{n,j} &= \left(\frac{2}{b}\right)^{j} \frac{n+1-j + \frac{2\beta}{b}(n-j)}{n+1}, \\
        -m_{12}t^{-1}_{1,j} + m_{11}t^{-1}_{n,j} &= \left(\frac{2}{b}\right)^{n-1-j} \frac{j + \frac{2\beta}{b}(j-1)}{n+1},
    \end{align*}
    which in turn gives
    \begin{align*}
        t^{-1}_{i,1} (m_{11} t^{-1}_{1,j} - m_{12} t^{-1}_{n,j}) + t^{-1}_{i,n} (-m_{12}t^{-1}_{1,j} + m_{11}t^{-1}_{n,j}) = \left(\frac{2}{b}\right)^{i+j} \left[ \left(1 + \frac{2\beta}{b}\right) \left( 1 - \frac{i+j}{n+1} + \frac{2ij}{(n+1)^2}\right) - \frac{\frac{2\beta}{b}}{n+1}\right]
    \end{align*}
    Furthermore,
    \begin{align*}
        \Delta  = m^{2}_{11} - m^{2}_{12} =  \frac{\left( 1+ \frac{2\beta}{b}\right) \left(n+1 + \frac{2\beta}{b}(n-1)\right)}{n+1}.
    \end{align*}
    Substitution of the last two expressions above and simplification results in the formula stated in the lemma.
    
\end{proof}

\begin{lemma}[Singularity]
    Let $b = \pm 2$. The matrix $\wtT_n$ is singular if $\wtb = \text{sgn}(b)$ and $\wtb = \text{sgn}(b)\frac{(n-3)}{(n-1)}$. Moreover, for $|\wtb| \geq 1$, $\Delta \geq 0$. 
\end{lemma}

\begin{proof}
The singularity of $\wtT_{n}$ occurs if the determinant $\Delta = 0$ or $m_{11} = \pm m_{12}$, which after using~\eqref{eq:elementM}, is equivalent to
\begin{align*}
1 + \frac{2\beta}{b} \frac{n}{n+1} = \pm \beta \left(\frac{2}{b}\right)^n \frac{1}{n+1}.
\end{align*}
Rearranging terms, with $\beta = \wtb - b$, results in
\begin{align*}
  \wtb = b - \frac{b (n+1)}{2} \frac{1}{n \mp \left(\frac{2}{b}\right)^{n-1}}.
\end{align*}
For $b = 2$, the above equation  gives $\wtb = 1$ or $\wtb = (n-3)/(n-1)$. For $b = -2$, the above equation gives $\wtb = -1$ and $\wtb = -(n-3)/(n-1)$. These lead to the first part of the lemma.

The determinant $\Delta$ can be simplified to a quadratic equation in $\wtb$, for $b = \pm 2$:
   \begin{align*}
        \Delta = \frac{n-1}{n+1}\left(\wtb \mp 1 \right) \left(\wtb \mp \frac{n-3}{n-1} \right).
    \end{align*}
    The second part of the lemma is obtained by checking the sign of $\Delta$, noting that $(n-3)/(n-1)<1$.

\end{proof}

\begin{lemma}[General Trace]
    Let $\wtT_{n}$ be invertible with $b = \pm 2$ and $\beta = \wtb - b$. The trace of the inverse matrix $\wtT^{-1}_{n}$ is 
    \begin{align}
        \text{Tr}(\wtT^{-1}_{n}) = \frac{n(n+2)}{3b} -  \frac{2\beta n}{3\left(1+ \frac{2\beta}{b} \right)} + \frac{\beta n}{3\left(1 + \frac{2\beta}{b} \right) \left( n+1 + \frac{2\beta}{b}(n-1) \right)}.
    \end{align}
\end{lemma}
\begin{proof}
    We use the general formula for the trace of the inverse matrix~\eqref{eq:traceTt}. First, from~\eqref{tij_general} we have
    \begin{align*}
        t^{-1}_{i,i} = \left(\frac{2}{b}\right) \frac{i(n+1-i)}{n+1}, \quad t^{-1}_{i,1} = \left( \frac{2}{b}\right)^{i} \frac{n+1-i}{n+1}, \quad \text{and}\quad t^{-1}_{n,i} = \left(\frac{2}{b}\right)^{n+1-i}\frac{i}{n+1}.
    \end{align*}
    Using the partial sum $\sum_{i=1}^n i = \frac{1}{2}n(n+1)$ and $\sum_{i=1}^n i^2 = \frac{1}{6}n(n+1)(2n+1)$, one can show that
    \begin{align*}
        \sum_{i=1}^{n} t^{-1}_{i, i} = \frac{n(n+2)}{3b}, \quad
        \sum_{i=1}^{n} \left((t^{-1}_{i, 1})^2 + (t^{-1}_{n,i})^2 \right) = \frac{n(2n+1)}{3(n+1)}, \quad \text{and} \quad
        \sum_{i=1}^{n} t^{-1}_{i, 1} t^{-1}_{n,i} = \left(\frac{2}{b}\right)^{n+1} \frac{n(n+2)}{6(n+1)}.
    \end{align*}
    By using $m_{11}$ and $m_{12}$ given in the proof of Lemma~\ref{lemma:prem1}, we get
    \begin{align*}
        m_{11} \sum_{i=1}^{n} \left((t^{-1}_{i,1})^2 + (t^{-1}_{n,i})^2 \right) - 2 m_{12} \sum_{i=1}^{n} t^{-1}_{i,1} t^{-1}_{n,i} 
        = \frac{n}{3(n+1)} \left[2n+1 + \frac{4\beta}{b} (n-1)\right]
    \end{align*}
    Substitution of the above equation, the sum $\sum_{i=1}^n t^{-1}_{i,i}$, and the determinant formula $\Delta$ (see in the proof of Lemma~\ref{lemma:prem1}) into~\eqref{eq:traceTt} gives the trace formula stated in the lemma.
\end{proof}

Now, we will consider each case separately and derive their row sums and upper bounds for the norms of their inverse matrices.


\section{The case $b = 2$} \label{sec:b=2}


In this section, we shall derive bounds for the rowsum and for some norms of the inverse of $\wtT_n$, with $b = 2$. 

\subsection{Bounds for rowsums of the inverse}
To begin with, we note that in this case the formula~\eqref{gen_wtt} can be simplified into
\begin{align}
    \wtt_{i,j}^{-1} &= \frac{\left(j(1+\beta) - \beta\right) \left((n-i)(1+\beta) + 1\right)}{(1+\beta)(n+1+\beta(n-1))} = \left(j - \gamma \right) \frac{(n-i)(1-\wtb)-1}{(n-1)(1-\wtb)-2}, \quad i \ge j, \label{eq1_wtt}
\end{align}
with $\gamma := \frac{2 - \wtb}{1 - \wtb}$.

\begin{theorem}[Row Sums] \label{lem1_rowsum}

Let $b = 2$ and $\wtb$ be such that $T_n$ is nonsingular. The $i$-th row sum of $\wtT^{-1}_n$ is given by
\begin{align*}
   \text{rowsum}_i \wtT^{-1}_n = \sum_{j=1}^{n} \wtt^{-1}_{i, j} =  \frac{1}{2} i(n+1-i) - \frac{n}{2}\frac{\wtb-2}{\wtb-1}.
\end{align*}
\end{theorem}

\begin{proof}
Due to symmetry of $\wtT_n^{-1}$, 
\begin{align*}
  \sum_{j=1}^n \wtt^{-1}_{i,j} = \sum_{j=1}^i \wtt^{-1}_{i,j} + \sum_{j=i+1}^n \wtt^{-1}_{i,j} = \sum_{j=1}^i \wtt^{-1}_{i,j} + \sum_{j=i+1}^n \wtt^{-1}_{j,i}.
\end{align*}
Thus, using~\eqref{eq1_wtt},
\begin{align*}
 \sum_{j=1}^n \wtt^{-1}_{i,j} &= \sum_{j=1}^n \left(j - \gamma \right) \frac{(n-i)(1-\wtb)-1}{(n-1)(1-\wtb)-2} + \sum_{j=i+1}^n \left(i - \gamma \right) \frac{(n-j)(1-\wtb)-1}{(n-1)(1-\wtb)-2} \\
 &= \frac{(n-i)(1-\wtb)-1}{(n-1)(1-\wtb)-2} \left[ \frac{i(i+1)}{2} - \gamma i \right] + \frac{(i- \gamma)(n-i)}{(n-1)(1-\wtb)-2} \left[ \frac{(n-i-1)(1-\wtb)}{2} - 1 \right],
\end{align*}
after using the partial sum $\sum_{j=1}^n j = \frac{1}{2}n(n+1)$ and $\sum_{j=i+1}^n j = \frac{1}{2}(n(n+1)-i(i+1))$. Simplification of the right-hand expression leads to the lemma.

\end{proof}

\begin{theorem}[Bounds for Row Sums] \label{theorem_bounds_rowsum}
Let $b = 2$ and $\wtb$ be such that $\wtT_n$ is nonsingular. The following inequalities hold true:
$$
   \frac{n}{2(\wtb-1)} \le \text{rowsum}_i \wtT^{-1}_n \le \frac{(n+1)^2}{8} - \frac{n(\wtb-2)}{2(\wtb-1)}.
$$
\end{theorem}
\begin{proof}
The row sum in Theorem~\ref{lem1_rowsum} is a quadratic function of $i$. It attains its maximum at $i = (n+1)/2$ and its minimum at either $i = 1$ or $i = n$. Substituting these values of $i$ gives the upper and lower bounds, respectively.
\end{proof}

\subsection{Bounds for norms of the inverse}

In the context of finding bounds for the norm of the inverse matrix $\widetilde{T}_n^{-1}$, it's crucial to account for the change in signs of its entries when deriving these bounds. Initially, we'll establish a lower bound that can be improved specifically for the case when $\widetilde{b} < 1$. Following that, we will explore various cases to determine the respective upper bounds.

\begin{theorem}[Lower Bound] \label{LB_theorem}

Let $b = 2$ and $\wtb$ be such that $\wtT_n$ is invertible. Then 
\begin{equation} \label{LB_1}
L:= \max \left\{ \left|\frac{n(n-2)}{8} - \frac{n}{2(1-\widetilde{b})} \right|, \left| \frac{n}{2(1-\widetilde{b})} \right| \right\} \leq \|\widetilde{T}_n^{-1}\|_{\infty}.
\end{equation}
\end{theorem}
\begin{proof}

Because of the change of signs along a row of $\wtT_n^{-1}$, we have
\begin{align*}
\|\wtT^{-1}_{n}\|_{\infty} = \max_{i} \sum_{j=1}^{n} \left| \wtt^{-1}_{i,j}\right| \geq \max_{i} \left|\sum_{j=1}^{n} \wtt^{-1}_{i,j} \right| \geq \max \left\{\left|\max_{i}{\sum_{j=1}^{n} \wtt^{-1}_{i,j}}\right|, \left|\min_{i}{\sum_{j=1}^{n} \wtt^{-1}_{i,j}}\right|\right\},
\end{align*}
with the min part given by Theorem~\ref{lem1_rowsum}. For the max part, consider the rowsum of $\wtT^{-1}_{n}$ in Theorem~\ref{lem1_rowsum} as a function of $i$. For all integer $n$,
\begin{align*}
    \sum_{j=1}^{n}  \wtt^{-1}_{\frac{n}{2}, j} = \frac{n(n+2)}{8} - \frac{n}{2} \frac{\wtb - 2}{\wtb - 1} \leq  \frac{(n+1)^2}{8} - \frac{n}{2} \frac{\wtb - 2}{\wtb - 1} = \sum_{j=1}^{n}  \wtt^{-1}_{\frac{n+1}{2}, j},
\end{align*}
which gives the max part. Note that
\begin{align*}
    \frac{n(n+2)}{8} - \frac{n}{2} \frac{\wtb - 2}{\wtb - 1} = \frac{n(n-2)}{8} - \frac{n}{2(1-\wtb)}.
\end{align*}

\begin{theorem}[Upper Bound]
For $b = 2$ and $\wtb$ such that $\wtT_{n}$ is invertible, we have
 \begin{equation}
        \|\wtT^{-1}_{n}\|_{\infty} \leq 
        U := \begin{cases}
            \frac{(n+1)^2}{8} - \frac{n\gamma}{2} & \text{for } \wtb > 1,\\
            \frac{n(\gamma-1)}{2} & \text{for } \frac{n-2}{n-1} \leq \wtb < 1,\\
            \max\left\{P, Q \right\} & \text{for } \frac{n-3}{n-1} < \wtb < \frac{n-2}{n-1},\\
            \max\left\{-P, R \right\} & \text{for } 0 < \wtb < \frac{n-3}{n-1},\\
            \max\left\{S, T \right\} & \text{for } \wtb \leq 0
        \end{cases}
    \end{equation}
    where 
    \begin{align*}
        P &= \frac{n(1-\gamma)}{2} + \frac{(\gamma-1)^2 (\gamma + 1)}{2\gamma - n - 1}, \\
        Q &= \frac{n(\gamma-1)}{2} + \frac{\gamma}{2\gamma - n -1}\frac{(n+1)^2}{16} + \frac{1}{2}, \\
        R &= \frac{(n+1)^2}{8} - \gamma \left(\frac{n+1}{2} -\gamma \right), \\
        S &= \frac{(n+1)^2}{8} + \frac{4 - n(2-\wtb)}{2(1-\wtb)}, \\
        T &= \frac{1}{1-\wtb} \left( \frac{n}{2} + \frac{2}{(n-1)(1-\wtb) - 2} \right)
    \end{align*}
    with $\gamma = \frac{2 - \wtb}{1 - \wtb}$.
\end{theorem}

\end{proof}

\subsection{The subcase $\wtb  >  1$}

In this subcase, the matrix $\widetilde{T}_n$ is symmetric and positive definite. Therefore, $\widetilde{T}_n^{-1} > 0$ and hence $|\widetilde{t}_{i,j}^{-1}| = \widetilde{t}_{i,j}^{-1}$. Consequently, we have
\[
\sum_{j=1}^n |\wtt_{i,j}^{-1}| = \text{rowsum}_i \wtT^{-1}_n.
\]
By using Theorem~\ref{theorem_bounds_rowsum}, the norm of the inverse matrix can be bounded from above as follows:
\begin{align} \label{c1_norm}
\|\wtT^{-1}_n\|_{\infty} = \max_i \sum_{j=1}^n |\wtt_{i,j}^{-1}| = \max_i \text{rowsum}_i \wtT^{-1}_n \le \frac{(n+1)^2}{8} -  \frac{n(\wtb-2)}{2(\wtb-1)}.
\end{align}

\subsection{The subcase $\wtb \le 0$}

We first derive the following lemma, which characterizes the signs of elements of the inverse matrix $\wtT_n^{-1}$. 

\begin{lemma} \label{lemma:btmin}
Let $b = 2$ and $\wtb < 0$. Then
\begin{align}
\text{sgn}\left(\wtt_{ij}^{-1}\right) =
\begin{cases}
  +1,& \text{for } (i,j) \in \{2,\dots,n-1\}^2  \cup \{(1,n), (n,1)\} \\
  -1,&\text{otherwise}.
\end{cases} \notag
\end{align}
\end{lemma}

\begin{proof}
We start by rewriting~\eqref{eq1_wtt} as 
\begin{align}
\wtt_{i,j}^{-1} = \left(j - \gamma \right) \frac{(n-i)(1+|\wtb|)-1}{(n-1)(1+|\wtb|)-2}, \quad i \ge j
\end{align}
with $\gamma = \frac{2+|\wtb|}{1+|\wtb|} = 1 + \frac{1}{1+|\wtb|} \in (1,2)$. For $|\wtb| > 0$, $(n-1)(1+|\wtb|)-2 > (n-1)-2 = n-3 \ge 0$, for $n \ge 3$. Furthermore, $(n-i)(1+|\wtb|) - 1 > n - i - 1 \ge 0$ for $i = 1,\dots,n-1$ and $j-\gamma > j - 2 \ge 0$, 
for $j = 2,\dots,n$. Thus, for $(i,j) \in \{2,\dots,n-1\} \times \{2,\dots,n-1\}$, we have $\wtt^{-1}_{i,j} > 0$.

Considering the case when $j = 1$, we have $j - \gamma = - \frac{1}{1+|\wtb|} < 0$. Therefore, $\wtt_{i,1}^{-1} < 0$, for $i = 1,\dots, n-1$. For $i = n$, $(n-i)(1+|\wtb|) - 1 < 0$. Thus, with $j - \gamma < 0$, we get $\wtt^{-1}_{n,1} > 0$.

The proof is completed by using the fact that $\wtT_n^{-1}$ is symmetric and centrosymmetric; i.e., $\wtt_{i,j}^{-1} = \wtt_{j,i}^{-1}$ and $\wtt_{i,j}^{-1} = \wtt_{n-i+1,n-j+1}^{-1}$.
\end{proof}

Using the same line of proof as above, the following result can be obtained:
\begin{lemma} \label{lemma:bt0}
Let $b = 2$ and $\wtb = 0$. Then
\begin{align}
\text{sgn}\left(\wtt_{ij}^{-1}\right) =
\begin{cases}
  +1,& \text{for } (i,j) \in \{3,\dots,n-2\}^2  \cup \{(1,n), (n,1)\} \\
   0,& \text{for } (i,j) \in \{2,n-1\}\times\{1,\dots,n\}\cup \{1,\dots,n\}\times\{2,n-1\} - \{(1,2),(2,1),(n-1,n),(n,n-1)\} \\
  -1,&\text{otherwise}.
\end{cases} \notag
\end{align}
\end{lemma}

\begin{theorem}
Let $b = 2$ and $\wtb \le 0$. Then
\begin{align}
    \|\wtT_n^{-1}\|_{\infty} \le  \max \left\{ \frac{1}{8}(n+1)^2 + \frac{4 - (2+|\wtb|)n}{2(1+|\wtb|)}, \frac{1}{1+|\wtb|} \left( \frac{n}{2} + \frac{2}{(n-1)(1+|\wtb|) - 2} \right) \right\}. \notag
\end{align}
\end{theorem}

\begin{proof}
For $i = 2,\dots,n-1$,
\begin{align} 
  \sum_{j=1}^n |\wtt_{i,j}^{-1}| &= |\wtt^{-1}_{i,1}| + \sum_{j=2}^{n-1} |\wtt^{-1}_{i,j}| + |\wtt^{-1}_{i,n}| = |\wtt^{-1}_{i,1}| + \sum_{j=2}^{n-1} \wtt^{-1}_{i,j} + |\wtt^{-1}_{i,n}| = 2|\wtt^{-1}_{i,1}| + \sum_{j=1}^{n} \wtt^{-1}_{i,j} + 2|\wtt^{-1}_{i,n}|, \notag 
\end{align}
because of Lemmas~\ref{lemma:btmin} and~\ref{lemma:bt0}. With
\begin{align}
    \wtt_{i,1}^{-1} = -\frac{1}{1+|\wtb|} \cdot \frac{(n-i)(1+|\wtb|)-1}{(n-1)(1+|\wtb|)-2}, \quad
    \wtt_{i,n}^{-1} = \wtt_{n,i}^{-1} = - \frac{1}{1+|\wtb|} \cdot \frac{(i-1)(1+|\wtb|)-1}{(n-1)(1+|\wtb|)-2}, \notag
\end{align}
and Theorem~\ref{lem1_rowsum}, the above sum can be written as
\begin{align} 
  \sum_{j=1}^n |\wtt_{i,j}^{-1}| &= \frac{2}{1+|\wtb|} \left( \frac{(n-i)(1+|\wtb|)-1}{(n-1)(1+|\wtb|)-2} +  \frac{(i-1)(1+|\wtb|)-1}{(n-1)(1+|\wtb|)-2} \right) + \frac{1}{2} i(n+1-i) - \frac{n}{2}\frac{2+|\wtb|}{1+|\wtb|}. \notag
\end{align}
Treating the right-hand side as a continuous function of $i$ in the interval $[2,n-1]$, it can be shown that the maximum is attained at $i = (n+1)/2$. In this case,
\begin{align}
    \sum_{j=1}^n |\wtt_{\frac{n+1}{2},j}^{-1}| = \frac{1}{8}(n+1)^2 + \frac{4 - (2+|\wtb|)n}{2(1+|\wtb|)}, \quad j = 2,\dots,n-1. \notag
\end{align}

For $i = 1$, 
\begin{align}
    \sum_{j=1}^n |\wtt_{1,j}^{-1}| &= \sum_{j=1}^{n-1} |\wtt^{-1}_{1,j}| + |\wtt^{-1}_{1,n}| = - \sum_{j=1}^{n-1} \wtt^{-1}_{1,j} + \wtt^{-1}_{1,n} = - \sum_{j=1}^{n} \wtt^{-1}_{1,j} + 2 \wtt^{-1}_{n,1} \notag
\end{align}
due to symmetry and Lemmas~\ref{lemma:btmin} and~\ref{lemma:bt0}. By using
\begin{align}
    \wtt_{n,1}^{-1} = \frac{1}{(1+|\wtb|)[(n-1)(1+|\wtb|)-2]}, \notag
\end{align}
obtained from~\eqref{eq1_wtt}, and Theorem~\ref{lem1_rowsum}, the sum can be written as follows: 
\begin{align}
    \sum_{j=1}^n |\wtt_{1,j}^{-1}| &= \frac{n}{2(1+|\wtb|)} + \frac{\textcolor{blue}{2}}{(1+|\wtb|)[(n-1)(1+|\wtb|)-2]} = \frac{1}{1+|\wtb|} \left( \frac{n}{2} + \frac{\textcolor{blue}{2}}{(n-1)(1+|\wtb|) - 2} \right). \notag
\end{align}
Centrosymmetry of $\wtT^{-1}_n$ leads to $\sum_{j=1}^n |\wtt_{n,j}^{-1}| = \sum_{j=1}^n |\wtt_{1,j}^{-1}|$. Therefore,
$$
  \|\wtT^{-1}_n \|_{\infty} = \max_i \sum_{j=1}^n |\wtt^{-1}_{i,j}| = \max \left\{ \sum_{j=1}^n |\wtt^{-1}_{\frac{n+1}{2},j}|, \sum_{j=1}^n |\wtt^{-1}_{1,j}| \right\},
$$
as in the theorem.
\end{proof}

\begin{corollary}
 Let $b = 2$. For $n \ge 9$,
\begin{align}
    \|\wtT_n^{-1}\|_{\infty} \le  \frac{1}{8}(n+1)^2 + \frac{4 - (2+|\wtb|)n}{2(1+|\wtb|)}, \notag
\end{align}
for all $\wtb \le 0$, where equality holds when $n$ is odd.
\end{corollary}

\begin{proof}
For $n \ge 9$, $(n-1)(1+|\wtb|)-2 \ge 8(1+|\wtb|)-2 = 6 + 8|\wtb|$. Thus,
\begin{align}
    f_{d} &:=  \sum_{j=1}^n |\wtt^{-1}_{\frac{n+1}{2},j}| - \sum_{j=1}^n |\wtt^{-1}_{1,j}| = \frac{1}{8}(n+1)^2 + \frac{4 - (2+|\wtb|)n}{2(1+|\wtb|)} - \frac{1}{1+|\wtb|} \left( \frac{n}{2} + \frac{2}{(n-1)(1+|\wtb|) - 2} \right) \notag \\
    &\ge \frac{1}{8}(n+1)^2 + \frac{4 - (2+|\wtb|)n}{2(1+|\wtb|)} - \frac{1}{1+|\wtb|} \left( \frac{n}{2} + \frac{2}{8|\wtb| + 6} \right) \notag \\
    &= \frac{1}{8}(n+1)^2 - \frac{3+|\wtb|}{2(1+|\wtb|)} n + \frac{2(8|\wtb|+5)}{(1+|\wtb|)(8|\wtb| + 6)}.    \notag
\end{align}
For $n = 9$,
\begin{align}
    f_d(9,\wtb) = \frac{64|\wtb|^2 + 56|\wtb| + 4}{8|\wtb|^2 + 14|\wtb| + 6} > 0,\quad \text{for } \wtb \le 0.\notag
\end{align}
For $n > 9$,
\begin{align}
    f_d(n,\wtb) - f_d(9,\wtb) = \frac{(4|\wtb| + 3)(n-9)((n+7)|\wtb| + n - 1)}{4(1+|\wtb|)(8|\wtb|+6)} \ge \frac{ 3(n-9)(n-1)}{4(1+|\wtb|)(8|\wtb|+6)} > 0,\quad \text{for } \wtb \le 0.\notag
\end{align}
Consequently, $
   \max \left\{ \sum_{j=1}^n |\wtt^{-1}_{\frac{n+1}{2},j}|, \sum_{j=1}^n |\wtt^{-1}_{1,j}| \right\}
 =  \sum_{j=1}^n |\wtt^{-1}_{\frac{n+1}{2},j}|$.
\end{proof}

\subsection{The subcase $0 < \wtb < 1$}

In this case, the matrix $\wtT_n$ is no longer diagonally dominant and, for a given $n \ge 3$, is singular for some choice of $\wtb$. We first state the end result, summarized in the theorem below.

\begin{theorem}
 \label{theor_upper_bound_1}
    Let $b = 2$, $0 < \wtb < 1$, and $\gamma = (2-\wtb)/(1-\wtb)$ such that $\wtT_{n}$ is invertible. Then
    \begin{equation} \label{upper_bound_1}
        \|\wtT^{-1}_{n}\|_{\infty} \leq 
        \begin{cases}
            \displaystyle \frac{1}{2}(\gamma-1)n &\text{for $ \frac{n-2}{n-1} \le \wtb < 1$},\\
            \displaystyle \max\left\{\frac{n(1-\gamma)}{2} + \frac{(\gamma-1)^2 (\gamma +1)}{2\gamma - n- 1}, \frac{n(\gamma-1)}{2} + \frac{\gamma}{2\gamma - n -1}\frac{(n+1)^2}{16} + \frac{1}{2}\right\} &\text{for $\frac{n-3}{n-1} < \wtb < \frac{n-2}{n-1}$}, \\
            \displaystyle \max\left\{\frac{n(\gamma-1)}{2} + \frac{(\gamma-1)^2(\gamma+1)}{n+1-2\gamma}, \frac{(n+1)^2}{8} - \gamma \left( \frac{n+1}{2} - \gamma \right)\right\} &\text{for $0 < \wtb < \frac{n-3}{n-1}$}.
        \end{cases} \notag
    \end{equation}   
\end{theorem}

\begin{proof}
We shall show the proof for Case 1 and 2. For Case 3, the proof is presented in Appendix~\ref{App:2.2}.
\newline
\noindent $\bullet$ Case 1: $\frac{n-2}{n-1} \le \wtb < 1$.
The matrix $\wtT_n$ in this case is nonsingular. Consider the elements of the inverse matrix given by the formula~\eqref{eq1_wtt}. It can be shown $(n-i)(1-\wtb)-1 \le 0$ and $(n-1)(1-\wtb) - 2 < 0$, for $i = 1,\dots,n$. Furthermore, $j - \gamma = j - 1 - \frac{1}{1-\wtb} \leq n-1- \frac{1}{1-\wtb} = ((n-1)(1-\wtb)-1)/(1-\wtb)\leq 0$, for all $j = 1,\dots,n$. Thus, $\wtt_{i,j}^{-1} < 0$.

We then have
\begin{align*}
    \sum_{j=1}^n \left|\wtt_{i,j}^{-1}\right|&= \sum_{j=1}^i \left|\wtt_{i,j}^{-1}\right| + \sum_{j=i+1}^n \left|\wtt_{i,j}^{-1}\right| = \sum_{j=1}^i \left|\wtt_{i,j}^{-1}\right| + \sum_{j=i+1}^n \left|\wtt_{j,i}^{-1}\right| \\
      &= \sum_{j=1}^{i} \left| (j-\gamma) \frac{(n-i)(1-\wtb)-1}{(n-1)(1-\wtb)-2} \right| + \sum_{j=i+1}^{n} \left| (i-\gamma) \frac{(n-j)(1-\wtb)-1}{(n-1)(1-\wtb)-2} \right| \\
      &= - \left( \sum_{j=1}^i (j-\gamma) \frac{(n-i)(1-\wtb)-1}{(n-1)(1-\wtb)-2}  + \sum_{j=i+1}^n (i-\gamma) \frac{(n-j)(1-\wtb)-1}{(n-1)(1-\wtb)-2} \right) \\
      &= \frac{1}{2} i(i-(n+1)) + \frac{n}{2}\frac{\wtb-2}{\wtb-1},
\end{align*}
due to Theorem~\ref{lem1_rowsum}. The right-hand side term attains maximum at $i = 1$ or $i = n$, resulting in the first part of the theorem.
\newline
\noindent $\bullet$ Case 2: $\frac{n-3}{n-1} < \wtb < \frac{n-2}{n-1}$. We shall use the following formula for the absolute value of the elements of the inverse matrix, obtained directly from~\eqref{eq1_wtt}:
\begin{align}
  \left| \wtt^{-1}_{i, j}\right| &= \begin{cases}
            \frac{\left| j - \gamma\right|}{\left| 1 - \gamma \right|}  \left| \wtt^{-1}_{i, 1} \right|, &\text{for } i \ge j, \\
            \left|(n-j)(1-\wtb) -1 \right| \left| \wtt^{-1}_{i,n} \right|, &\text{for } i < j.
        \end{cases} \notag
\end{align}
Therefore,
\begin{align} \label{eq1_abs_sum}
    \sum_{j=1}^{n} \left| \wtt^{-1}_{i, j}\right| 
            = \sum_{j=1}^{i}\left| \wtt^{-1}_{i, j}\right| + \sum_{j=i+1}^{n} \left| \wtt^{-1}_{i, j}\right| = \frac{\left| \wtt^{-1}_{i, 1}\right|}{\left| 1- \gamma \right|} \sum_{j=1}^{i} \left|j-\gamma\right| + \left|\wtt^{-1}_{i,n} \right| \sum_{j=i+1}^{n} \left|(n-j)(1-\wtb)-1 \right|.
\end{align}

Consider the convex combination $\wtb = (1-\alpha) \frac{n-3}{n-1} + \alpha \frac{n-2}{n-1}$, with $\alpha \in (0,1)$. Then $(n-1)(1-\wtb)-2 = -\alpha <0$. The term $(n-i)(1-\wtb) - 1$ in $\wtt^{-1}_{i,j}$ is however not sign-definite. This term vanishes when $i = n - 1/(1-\wtb) = n - \gamma + 1 =: \widetilde{i}$. Let $i^* = \lfloor \widetilde{i} \rfloor$. Then, for $1 \le i \le i^*$, we have $(n-i)(1-\wtb) - 1 \ge 0$ and for $i^* < i \le n$ we have $(n-i)(1-\wtb) - 1 < 0$. Using the aforementioned convex combination, we obtain the following facts:
\begin{enumerate}
    \item[(i)] $\widetilde{i} \in \left(1,\frac{n+1}{2}\right)$. Thus, $i^* \in \left[1,\lfloor \frac{n+1}{2} \rfloor \right]$;
    \item[(ii)] $\gamma \in \left( \frac{n+1}{2}, n \right)$. Thus, $\lfloor \gamma \rfloor \in \left[ \lfloor \frac{n+1}{2} \rfloor, n-1\right] $.
\end{enumerate}
Consequently, $i^* \le \lfloor \frac{n+1}{2} \rfloor \le \frac{n+1}{2} < \gamma$, and $i^* \le \lfloor \gamma \rfloor$. 

\paragraph{Subcase 2.1:} $1 \le i < i^*$.
From Fact (i),
$$
 1- \gamma < \dots < (i^*-1) - \gamma < i^* - \gamma < 0.
$$
For first term in the sum~\eqref{eq1_abs_sum}, since $i < i^* < \gamma$, then $j - \gamma < 0$.
Hence,
$$
  \sum_{j=1}^{i} |j-\gamma| = \sum_{j=1}^i (\gamma - j) = \gamma i - \frac{i(i+1)}{2}.
$$
For the second term in the sum~\eqref{eq1_abs_sum}, we have
\begin{align}
    \sum_{j=i+1}^{n} \left|(n-j)(1-\wtb)-1 \right| \le 1 + \gamma + (n-i) \left[ \frac{1}{2} (n-i-1)(1-\wtb) - 1\right]. \label{eq:subcase2.1second}
\end{align}
Detail of the derivation is presented in Appendix~\ref{App:2.1}. Thus, with
\begin{align*}
\left| \wtt^{-1}_{i, 1} \right| &= \left| (1-\gamma) \frac{(n-i)(1-\wtb)-1}{(n-1)(1-\wtb)-2} \right| = (\gamma - 1)\frac{(n-i)(1-\wtb)-1}{2 - (n-1)(1-\wtb)},\\
\left| \wtt^{-1}_{i,n} \right| &= \left| \wtt^{-1}_{n, i} \right| = \left| \frac{i-\gamma}{(n-1)(1-\wtb)-2} \right| = \frac{\gamma - i}{2 - (n-1)(1-\wtb)},
\end{align*}
and Lemma~\ref{lem1_rowsum} and its proof, it can be shown that
\begin{align}
    \sum_{j=1}^n |\wtt_{i,j}^{-1} |  &\le \sum_{j = 1}^{-1} \wtt_{i,j}^{-1} + \frac{(\gamma + 1)(i-\gamma)}{(n-1)(1-\wtb)-2} = \frac{1}{2} i(n+1-i) - \frac{n}{2}\gamma + \frac{(\gamma + 1)(i-\gamma)}{(n-1)(1-\wtb)-2} =: P(i).  \label{eq:Pi}
\end{align}
For $i \ge 1$,
\begin{align}
    P'(i) \leq -1 + \frac{n+1}{2} + \frac{\gamma+1}{(n-1)(1-\wtb)-2} = \frac{[(n-1)(1-\wtb) - 1]^2 +4(1-\wtb)+1}{2(1-\wtb)[(n-1)(1-\wtb)-2]} < 0. \notag
\end{align}
Thus, $P(i)$ is a decreasing function, with maximum attained at $i = 1$. The following bound then results:
\begin{align*}
\max_{1\leq i \leq i^*} \sum_{j=1}^{n}\left| \wtt^{-1}_{i, j} \right|  \leq P: = \max_{i} P_{i} = P_{1}  = \frac{n(1-\gamma)}{2} - \frac{\gamma^2-1}{(n-1)(1-\wtb)-2} = \frac{n(1-\gamma)}{2} + \frac{(\gamma-1)^2 (\gamma +1)}{2\gamma - n- 1}.
\end{align*}

\paragraph{Subcase 2.2:}
$i = i^*$. Performing similar calculations as in Subcase 2.1 and substituting to \eqref{eq1_abs_sum}, we get 
\begin{align*}
    \sum_{j=1}^{n}\left| \widetilde{t}^{-1}_{i^*, j} \right| &= \frac{(n-i^*)(1-\widetilde{b})-1}{(n-1)(1-\widetilde{b})-2} \left[\frac{i^*(i^* + 1)}{2} - \gamma i^* \right] - \frac{(i^* - \gamma)(n-i^*)}{(n-1)(1-\widetilde{b})-2} \left[ \frac{(n-i^*-1)(1-\widetilde{b})}{2} - 1\right] \\
    &= \frac{(n+1-i^*-\gamma)i^*(i^* - 2\gamma + 1)}{2(n+1-2\gamma)} - \frac{(i^* - \gamma)(n-i^*)(n+1-i^*-2\gamma)}{2(n+1-2\gamma)}\\
    &= \frac{1}{2}\left( \frac{i^* (n - i^*)(n+1-2i^*)}{2\gamma - n - 1} + \gamma(n - 2i^*) + i^*\right).
\end{align*}
Let us denote 
$$G_{i^*} := \frac{i^* (n - i^*)(n+1-2i^*)}{2\gamma - n - 1} - 2 \gamma i^* + i^*.$$
Now, we seek an upper bound for $T_{i^*}$ using $n - i^* < \gamma$
\begin{align*}
    G_{i^*} &< \frac{\gamma i^* (n+1-2i^*)}{2\gamma - n - 1} - 2\gamma i^* + i^* \\
    &\leq \frac{\gamma}{2\gamma - n - 1} \max_{i^*}\{i^*(n+1-2i^*)\} + \max_{i^*}\{i^* - 2\gamma i^*\}\\
    &\leq \frac{\gamma}{2\gamma - n - 1} \frac{(n+1)^2}{8} + 1 - 2\gamma\\
    &\leq \frac{\gamma}{2\gamma - n - 1} \frac{(n+1)^2}{8} + 1 - n.
\end{align*}
Therefore,
\begin{align*}
    \sum_{j=1}^{n}\left| \widetilde{t}^{-1}_{i^*, j} \right| = \frac{1}{2}\left(G_{i^*} + \gamma n \right) < Q: = \frac{n(\gamma-1)}{2} + \frac{\gamma}{2\gamma - n -1}\frac{(n+1)^2}{16} + \frac{1}{2}.
\end{align*}

\paragraph{Subcase 2.3:} $i^* < i \le \lfloor \frac{n+1}{2} \rfloor$.
In this case, the first sum in~\eqref{eq1_abs_sum} can be expressed as:
\begin{align*}
    \sum_{j=1}^n \left| j-\gamma \right| = \sum_{j=1}^n (\gamma - j) = \gamma i - \frac{1}{2}i(i+1).
\end{align*}
For the second sum, because $j \ge i+1 \ge i^*$,
\begin{align*}
    \sum_{j=i+1}^n \left|(n-j)(1-\wtb)-1\right| = \sum_{j=i+1}^n  (1- (n-j)(1-\wtb)) = \frac{1}{2}(n-i)(2-(1-\wtb)(n-i-1)).
\end{align*}
Using
\begin{align*}
    \left| \wtt_{i,1}^{-1} \right| = (\gamma-1) \frac{(n-i)(1-\wtb)-1}{(n-1)(1-\wtb)-2}, \quad |1-\gamma| = \gamma - 1, \quad \left| \wtt_{i,n}^{-1} \right| = \frac{(i-\gamma)(1-\wtb)-1}{(n-1)(1-\wtb)-2},
\end{align*}
the sum~\eqref{eq1_abs_sum} becomes
\begin{align*}
    \sum_{j=1}^{n} \left| \wtt^{-1}_{i, j}\right| &= -\frac{(n-i)(1-\wtb)-1}{(n-1)(1-\wtb)-2}\left(\frac{1}{2}i(i+1) - \gamma i \right) - \frac{(i-\gamma)(n-i)}{(n-1)(1-\wtb)-2} \left[1 - \frac{(n-i-1)(1-\wtb)}{2} \right] \\
    &= -\sum_{j=1}^n \wtt^{-1}_{i, j} = -\left( \frac{1}{2} i(n+1-i) - \frac{1}{2} \gamma n\right) =: R_i.
\end{align*}
With $R'_i = i - \frac{n+1}{2} < 0$, for $i \in (i^*, \lfloor \frac{n+1}{2} \rfloor]$, $R_i$ is a decreasing function of $i$, which attains minimum at $i = \frac{n+1}{2}$. Since $R_i$ is a quadratic function, $R_1 = \frac{n}{2}(\gamma-1) > R_i$, for all $i \in (i^*, \lfloor \frac{n+1}{2} \rfloor]$.

We now compare $P$ and $R_1$:
\begin{align*}
  P - R_1 &= \frac{n(1-\gamma)}{2} + \frac{(\gamma-1)^2 (\gamma +1)}{2\gamma - n- 1} - \frac{n}{2}(\gamma-1) = n(1-\gamma) +  \frac{(\gamma-1)^2 (\gamma +1)}{2\gamma - n- 1} \\
  &= \frac{(\gamma-1)^2 (\gamma + 1 - 2n) + (\gamma-1) n (n-1)}{2\gamma - n - 1} > 0
\end{align*}
for $n \ge 3$ and $\gamma \in (\frac{n+1}{2},n)$. Thus,
$\displaystyle \max_{1 \le i \le n} \sum_{j=1}^{n} \left| \wtt^{-1}_{i, j}\right| = \max \left\{P,Q, R_1\right\} = \max \left\{P, Q\right\}.$

\end{proof}

\section{The case $b = -2$} \label{sec:b=-2}

In this section, we shall consider the case where $b = -2$. Firstly, in this case, using~\eqref{tij_general}, the elements of the inverse Toeplitz matrix $T_n^{-1}$ are given by the formula
\begin{equation}
    t^{-1}_{i, j} = \begin{cases}
         (-1)^{i+1-j}  \frac{j(n+1-i)}{n+1} & \text{for } i \geq j\\
        t^{-1}_{j, i} & \text{for } i < j.
    \end{cases} \label{eq:Tijelement}
\end{equation}
Note that the rational term on the right-hand side of the above formula is the inverse elements when $b = 2$. Denote the latter by $b_{+}$ and let $T_{n,+}$ be the corresponding Toeplitz matrix. Then, for $b = -2 = -b_+$ and $i \ge j$, we can write the formula~\eqref{eq:Tijelement} as
\begin{align*}
    t^{-1}_{i, j} = (-1)^{i+1-j} t^{-1}_{i, j, +}.
\end{align*}

Let $\wtb = -\wtb_+$, where $\wtb_+ \in \mathbb{R}$ is the $(1,1)$- and $(n,n)$-element of the near Toeplitz matrix $\wtT_{n,+}$. 
The elements of $\wtT_n^{-1}$ with $b = -2$ is related to the elements of $\wtT_{n,+}^{-1}$ in a way as stated in the following lemma:
\begin{lemma} \label{lem:3_1}
    Let $b = -2 = -b_+$ and $\wtb = - \wtb_+$, with $\wtb_+$ be such that $\wtT_{n,+}$ is invertible. For $i \ge j$, the following holds:
\begin{align*}
    \widetilde{t}^{-1}_{i,j} = (-1)^{i+1-j}\widetilde{t}^{-1}_{i, j, +}
\end{align*}
\end{lemma}
\begin{proof}
Note that $\beta = \wtb - b = -\wtb_+ - (-b_+) = -(\wtb_+ - b_+) = -\beta_+$. Hence $\beta/b =  \beta_+/b_+$. Substitution of these relations to~\eqref{gen_wtt} in Lemma~\ref{lemma:prem1} leads to the lemma.
\end{proof}

\begin{theorem}[Row Sums] \label{rowsum_b2}
    Let $\wtb$ be such that, for $b = -2$, $\wtT_{n}$ is invertible. Then 
    \begin{align*}
        \text{rowsum}_{i} \wtT^{-1}_{n} = \sum_{j=1}^{n} \wtt^{-1}_{i, j} = 
        \begin{cases}
            \frac{(-1)^i - 1}{4} + \frac{(-1)^{i} (\beta n - i(\beta + 1))}{2(n+1-\beta(n-1))} \;&\text{for}\;n \mid 2,\\
            \frac{(-1)^i - 1}{4} + \frac{(-1)^{i}\beta}{2(1-\beta)} \;&\text{for}\;n \nmid 2,
        \end{cases}
    \end{align*}
    where $\beta = \wtb + 2$. 
\end{theorem}
\begin{proof}
    Our proof is based on the rowsum formula~\eqref{eq:rowsumTt}. First, we compute the rowsum of $T^{-1}_{n}$. For this, we have    \begin{align}\label{cb_1}
        \sum_{j=1}^{n} t^{-1}_{i, j} &= \sum_{j=1}^{i} t^{-1}_{i, j} + \sum_{j=i+1}^{n}t^{-1}_{j, i} \notag \\
        &= \sum_{j=1}^{i} (-1)^{i+1-j} t^{-1}_{i,j, +} + \sum_{j=i+1}^{n} (-1)^{j+1-i} t^{-1}_{j, i, +} \notag \\
        &= (-1)^{i+1} \left[ \frac{n+1-i}{n+1}\sum_{j=1}^{i} (-1)^{j} j  + \frac{i}{n+1} \sum_{j=i+1}^{n} (-1)^{j} (n+1-j)\right].
    \end{align}
    where we have used Lemma~\ref{lem:3_1}. We note that 
    \begin{align*}
        \sum_{j=1}^{i} (-1)^{j} j &= (-1)^{1}1 + (-1)^{2}2 + \ldots + (-1)^{i} i = \frac{(2i+1)(-1)^{i} - 1}{4}
    \end{align*}
    and 
    \begin{align*}
        \sum_{j=i+1}^{n} (-1)^{j} (n+1-j) = -\frac{(n-i)(-1)^{i}}{2} + \frac{(-1)^{n} - (-1)^{i}}{4}.
    \end{align*}
    Substitution of the above sums into \eqref{cb_1} and simplification results in 
    \begin{align*}
        \sum_{j=1}^{n} t^{-1}_{i, j} = \begin{cases}
            \frac{(-1)^{i} - 1}{4} - \frac{(-1)^{i} i}{2(n+1)} \;&\text{for}\;n \mid 2,\\
            \frac{(-1)^{i} - 1}{4} \; &\text{for}\;n \nmid 2. 
        \end{cases}
    \end{align*}
    If $i = n$, the sum becomes
    \begin{align*}
        \sum_{j=1}^{n} t^{-1}_{n, j} = \begin{cases}
            -\frac{n}{2(n+1)} &\text{for $n \mid 2$},\\
            -\frac{1}{2} &\text{for $n \nmid 2$}.
        \end{cases}
    \end{align*}
    Next,
    \begin{align*}
        t^{-1}_{i, 1} + t^{-1}_{i, n} &= (-1)^{i} \frac{n+1-i}{n+1} + (-1)^{n+1-i} \frac{i}{n+1} = 
        \begin{cases}
            (-1)^{i}(1 - \frac{2i}{n+1}) \; &\text{for}\; n \mid 2,\\
            (-1)^{i}  \; &\text{for}\; n \nmid 2.
        \end{cases}
    \end{align*}
    Finally, by evaluating the last formula for $i = n$, we have
    \begin{align*}
        m_{11} + m_{12} &= 1 + \beta(t^{-1}_{n, 1} + t^{-1}_{n, n}) = \begin{cases}
            1 -\frac{\beta(n-1)}{n+1} &\text{for $n \mid 2$},\\
            1-\beta &\text{for $n \nmid 2$}.
        \end{cases}
    \end{align*}
    Substitution of all above relations into \eqref{eq:rowsumTt} gives us the desired result. 
\end{proof}

\begin{theorem}[Upper Bound]
\label{theor_upper_bound_2}
For $b = -2$ and $\wtb$ such that $\wtT_{n}$ is invertible, we have
 \begin{equation}
        \|\wtT^{-1}_{n}\|_{\infty} \leq 
        U := \begin{cases}
            \frac{(n+1)^2}{8} - \frac{n\gamma_{+}}{2} & \text{for } \wtb < - 1,\\
            \frac{n(\gamma_{+}-1)}{2} & \text{for }  - 1 < \wtb \leq -\frac{n-2}{n-1},\\
            \max\left\{P, Q \right\} & \text{for } - \frac{n-2}{n-1}  < \wtb < - \frac{n-3}{n-1},\\
            \max\left\{-P, R \right\} & \text{for } -\frac{n-3}{n-1} < \wtb < 0,\\
            \max\left\{S, T \right\} & \text{for } \wtb \geq 0
        \end{cases}
    \end{equation}
    where 
    \begin{align*}
        P &= \frac{n(1-\gamma_{+})}{2} + \frac{(\gamma_{+}-1)^2 (\gamma_{+} + 1)}{2\gamma_{+} - n - 1}, \\
        Q &= \frac{n(\gamma_{+}-1)}{2} + \frac{\gamma_{+}}{2\gamma_{+} - n -1}\frac{(n+1)^2}{16} + \frac{1}{2}, \\
        R &= \frac{(n+1)^2}{8} - \gamma_{+} \left(\frac{n+1}{2} -\gamma_{+} \right), \\
        S &= \frac{(n+1)^2}{8} + \frac{4 - n(2+\wtb)}{2(1+\wtb)}, \\
        T &= \frac{1}{1+\wtb} \left( \frac{n}{2} + \frac{2}{(n-1)(1+\wtb) - 2} \right)
    \end{align*}
    with $\gamma_{+} = \frac{2 + \wtb}{1 + \wtb}$.
\end{theorem}

\begin{proof}
From Lemma~\ref{lem:3_1}, $|\wtt_{i,j}^{-1}| = |\wtt_{i,j,+}^{-1}|$. Thus,
$\|\wtT^{-1}_{n}\|_{\infty} = \max_{i} \sum_{j=1}^{n} |\wtt^{-1}_{i, j}| = \max_{i} \sum_{j=1}^{n} |\wtt^{-1}_{i, j, +}| = \|\wtT^{-1}_{n, +}\|_{\infty}$.
\end{proof}

\section{Numerical Experiments} \label{sec:num}

In this section, we analyze upper and lower bounds and compare the maximum convergence rates for the fixed-point iteration method in Fisher's problem.

The lower bound proposed in Theorem \ref{LB_theorem} demonstrates effectiveness when $b = 2$ and $ \frac{n-2}{n-1} \leq \widetilde{b} < 1$. Outside these conditions, there is room for optimizing the bound.
\begin{figure}[htbp]
    \centering
    \begin{subfigure}{0.45\textwidth}
        \centering
        \includegraphics[width=\textwidth]{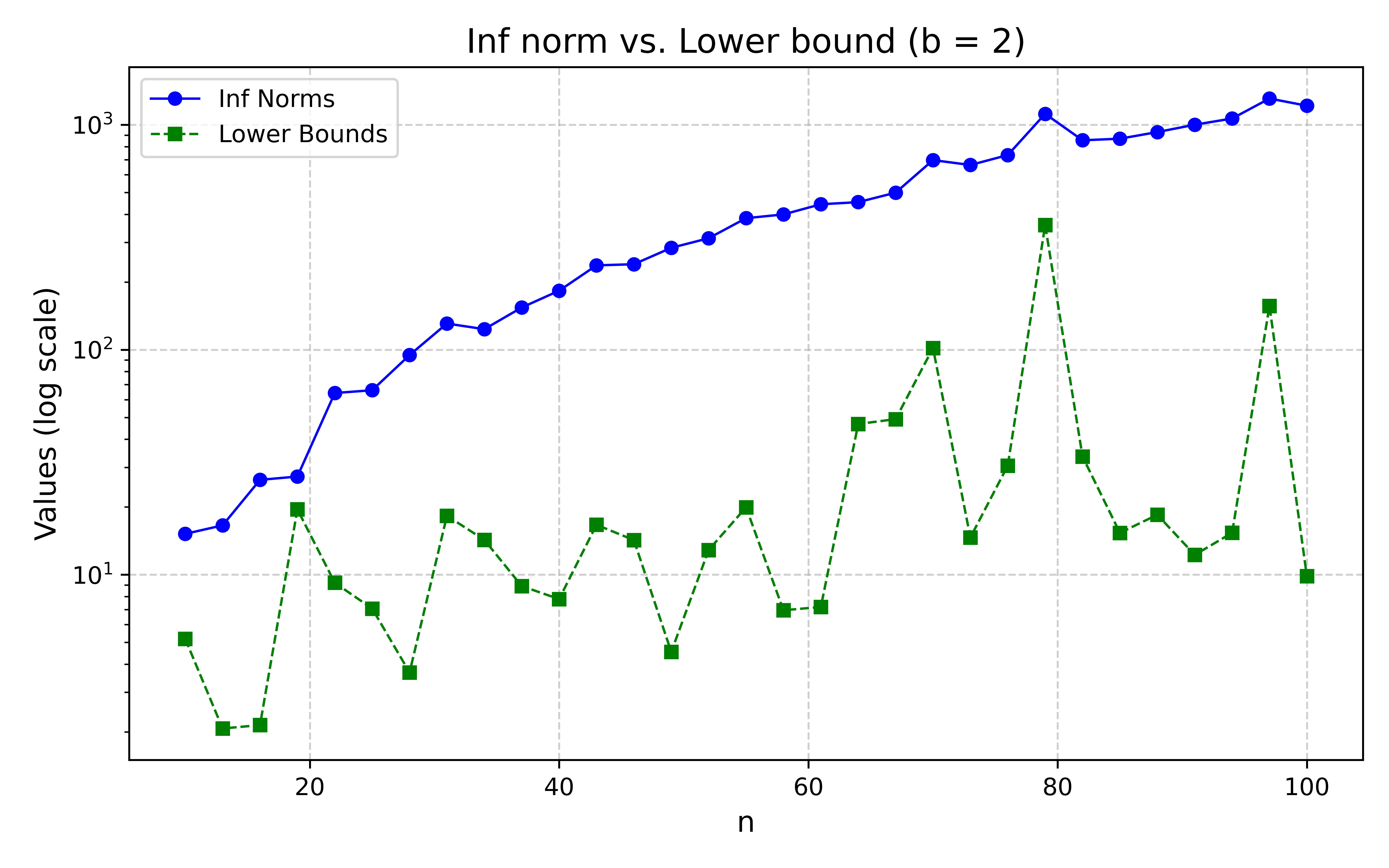}
        \caption{Log comparisons between $\|\wtT^{-1}_{n}\|_{\infty}$ and lower bounds.}
        \label{fig1}
    \end{subfigure}
    \hspace{\fill} 
    \begin{subfigure}{0.45\textwidth}
        \centering
        \includegraphics[width=\textwidth]{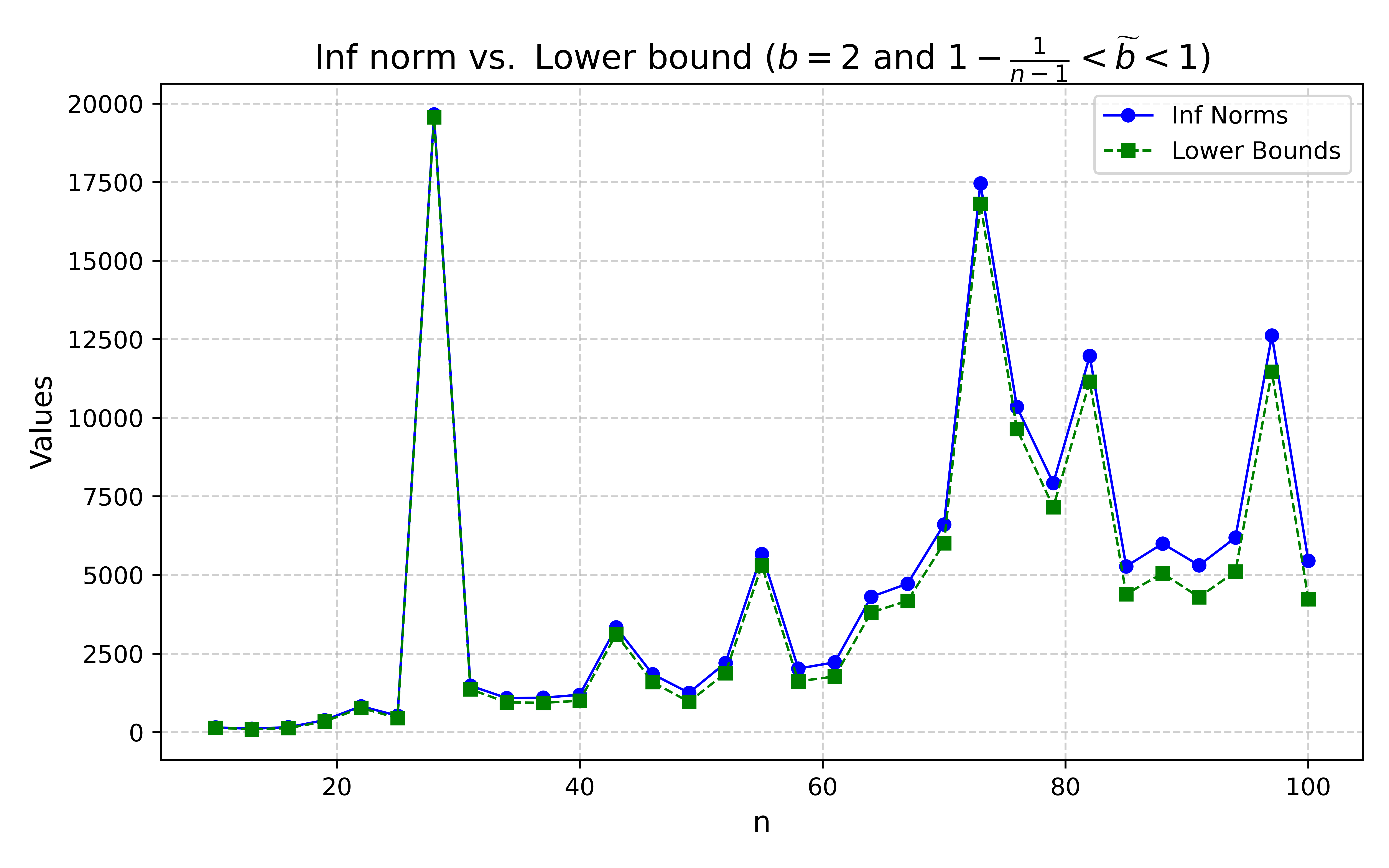}
        \caption{$\|\wtT^{-1}_{n}\|_{\infty}$ and lower bounds when $b = 2$ and $\frac{n-2}{n-1} \leq \wtb < 1$.}
        \label{fig11} 
    \end{subfigure}
    \caption{Comparison of $\|\wtT^{-1}_{n}\|_{\infty}$ vs. Lower bounds when $b = 2$.}
    \label{fig:combined1}
\end{figure}
\FloatBarrier
Following this, we proceed to illustrate upper bound examinations. Generally, we observe tight upper bounds when $b = 2$ with $\wtb \geq 1$ and $b = -2$ with $\wtb \leq -1$. In other cases, there is potential for improvement in alternate upper bounds. 
\begin{figure}[H]
    \centering
    \begin{subfigure}[b]{0.4\textwidth}
        \centering
        \label{tab1}
        \begin{tabular}{cccc}
            \toprule
            $n$ & $\widetilde{b}$ & $\|\wtT^{-1}_{n}\|_{\infty}$ & $\text{upper\_bound}$ \\
            \midrule
            10 & 5.93 & 11.014 & 11.139 \\
            13 & -2.28 &  16.627 & 24.500 \\
            16 & -3.46 &  26.656 & 36.125 \\
            19 & 3.03  & 45.188 &  45.188 \\
            22 & 6.39  & 57.041 & 57.166 \\
            25 & 0.46  & 54.197 & 55.506 \\
            28 & 11.62 & 92.319 & 92.444 \\
            \bottomrule
        \end{tabular}
        \caption{Comparison of $\|\wtT^{-1}_{n}\|_{\infty}$ and upper bounds.}
    \end{subfigure}%
    \hfill
    \begin{subfigure}[b]{0.45\textwidth}
        \centering
        \includegraphics[width=\textwidth]{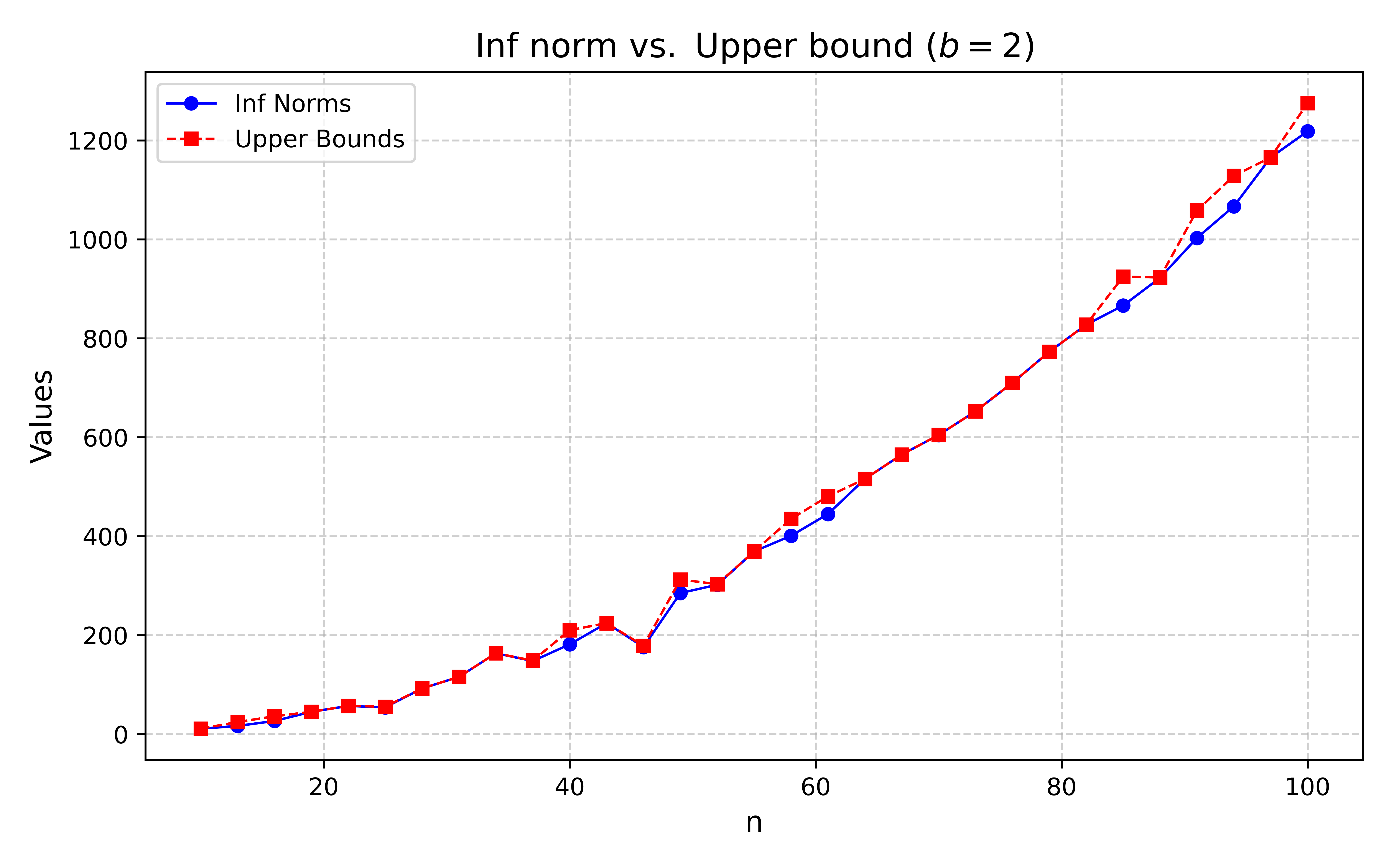}
        \caption{$\|\wtT^{-1}_{n}\|_{\infty}$ and upper bounds.}
        \label{fig2}
    \end{subfigure}
    \caption{Comparison of $\|\wtT^{-1}_{n}\|_{\infty}$ vs. Upper Bounds when $b = 2$.}
\end{figure}

\begin{figure}[H]
    \centering
    \begin{subfigure}[b]{0.4\textwidth} 
        \centering
        \label{tab5}
        \begin{tabular}{cccc}
            \toprule
            $n$ & $\widetilde{b}$ & $\|\widetilde{T}_{n}^{-1}\|_{\infty}$ & \text{upper\_bound} \\
            \midrule
            10 &  1.93 &  8.976 & 15.125 \\
            13 & -6.28 &  19.232 & 19.232 \\
            16 & -7.46 &  29.238 &  29.363 \\
            19 & -0.97 &  361.281 &  361.281 \\
            22 & 2.39  &  52.345 &  66.125 \\
            25 & -3.54 &  76.925 &  76.925 \\
            28 & 7.62 & 89.607 & 105.125 \\
            \bottomrule
        \end{tabular}
        \caption{Comparison of $\|\wtT^{-1}_{n}\|_{\infty}$ and upper bounds.}
    \end{subfigure}%
    \hfill
    \begin{subfigure}[b]{0.45\textwidth} 
        \centering
        \includegraphics[width=\textwidth]{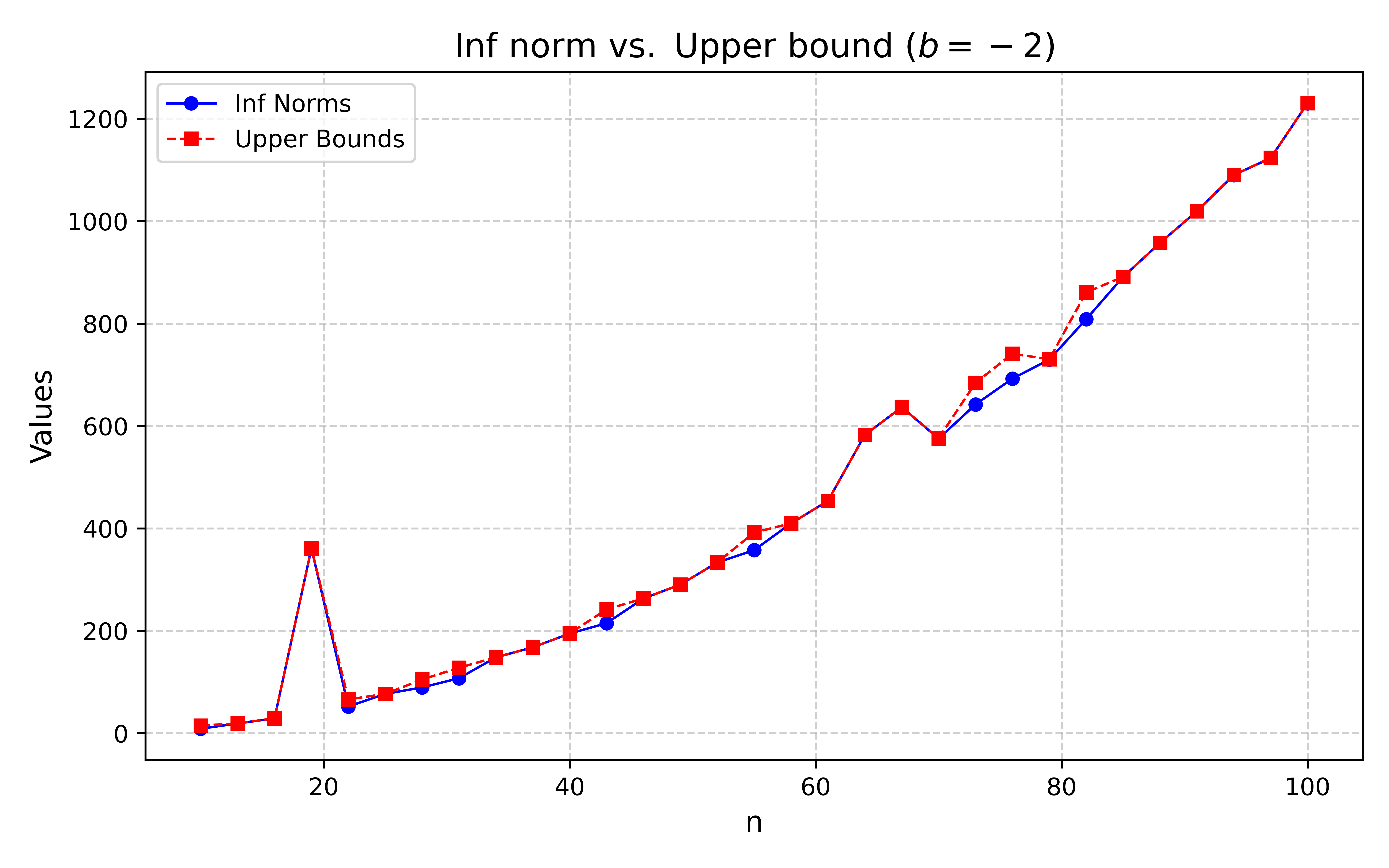} 
        \caption{$\|\widetilde{T}_{n}^{-1}\|_{\infty}$ and upper bounds.}
        \label{fig7}
    \end{subfigure}
    \caption{Comparison of $\|\widetilde{T}_{n}^{-1}\|_{\infty}$ vs. Upper bounds when $b = -2$.}
\end{figure}
\FloatBarrier
\begin{table}[H]
    \centering
    \caption{Observed maximum convergence rate for Fisher's problem with $n = 50$, $b = \wtb = 2.0$, and $L = 0.5$. The expected rates are determined using upper bounds from Theorem \ref{theor_upper_bound_1} applied to the norms of the inverse matrix.}
    \label{table5}
    \begin{tabular}{cccc}
        \toprule
        $k$ & $\textbf{Iterations}$ & $\textbf{Numerical rate}$ & $\textbf{Expected rate}$ \\
        \midrule
        1/2 & 5 & 0.0132 & 0.0163 \\
        1 & 5 & 0.0264 & 0.0325 \\
        2 & 6 & 0.0527 & 0.065 \\
        4 & 8 & 0.1054 & 0.13 \\
        8 & 10 & 0.2109 & 0.2601 \\
        16 & 17 & 0.4218 & 0.5202 \\
        32 & 68 & 0.8436 & 1.0404 \\
        \bottomrule
    \end{tabular}
\end{table}

\begin{table}[H]
    \centering
    \caption{Observed maximum convergence rate for Fisher's problem with $n = 50$, $b = \wtb = -2.0$, and $L = 0.05$. The expected rates are determined using upper bounds from Theorem \ref{theor_upper_bound_2} applied to the norms of the inverse matrix.}
    \label{table6}
    \begin{tabular}{cccc}
        \toprule
        $k$ & $\textbf{Iterations}$ & $\textbf{Numerical rate}$ & $\textbf{Expected rate}$ \\
        \midrule
        1 & 3 & 0.0003 & 0.0003 \\
        3 & 3 & 0.0006 & 0.001 \\
        9 & 3 & 0.002 & 0.0029 \\
        27 & 4 & 0.0071 & 0.0088 \\
        81 & 4 & 0.0214 & 0.0263 \\
        243 & 6 & 0.0641 & 0.079 \\
        729 & 9 & 0.1922 & 0.237 \\
        \bottomrule
    \end{tabular}
\end{table}

\section{Conclusion}

The objective of this paper is to provide properties of inverse of $\wtT$ in (\ref{eq:matrixTn}) with $|b| = 2$, where Toeplitz part with weakly diagonally dominant. We have derived bounds for both dominant ($|\widetilde{b}| \geq 1$) and non-dominant ($|\widetilde{b}| < 1$) cases. Numerical examples for Fisher's problem is shown to verify the obtained theoretical results.

\FloatBarrier
\bibliographystyle{plain}
\bibliography{reference}

\newpage
\appendix

\section{Part of Proof of Theorem~\ref{theor_upper_bound_1}}

\subsection{Case 2: Inequality~\eqref{eq:subcase2.1second}} \label{App:2.1}
For $1 \le i < i^*$ and $j = i+1, ..., n$,
\begin{align} 
\sum_{j=i+1}^{n} \left| (n-j)(1-\wtb)-1 \right| 
    &= \sum_{j=i+1}^{i^*} \left[ (n-j)(1-\wtb)-1 \right] + \sum_{j=i^*+1}^{n} \left[ 1 - (n-j)(1-\wtb) \right] \notag \\
    &= (i^* - i)(n(1-\wtb)-1) - \frac{(i^*-i)(i^* + i + 1)(1-\wtb)}{2} \notag\\
    &+ (n-i^*)(1 - n(1-\wtb)) + \frac{(n-i^*)(n+i^*+1)(1-\wtb)}{2}\notag \\
    &= (n-2i^*+i) - (1-\wtb)(n-i^*)(n-i^*-1)  + \frac{(n-i)(n-i-1)(1-\wtb)}{2} \notag
\end{align}
Because $(n-i^*)(1 - \wtb) \geq 1$ and that
\begin{align*}
    (n - i^* - 1)(1 - \wtb) - 1 < 0 \Rightarrow (n - i^*)(1 - \wtb) < 2 - \wtb \quad \Rightarrow \quad n - i^* < \gamma  \quad \text{or} \quad i^* > n - \gamma,
\end{align*}
the following inequalities are deduced:
\begin{align}\label{eq:App1.1}
\sum_{j=i+1}^{n} \left| (n-j)(1-\wtb)-1 \right|  &\leq (n-2i^*) - (n-i^* - 1) + i + \frac{(n-i)(n-i-1)(1-\wtb)}{2} \notag \\
    &\leq 1 - i^* + i + \frac{(n-i)(n-i-1)(1-\wtb)}{2}  \notag\\
    &\leq 1 + \gamma - n + i + \frac{(n-i)(n-i-1)(1-\wtb)}{2} \notag\\
    &\leq 1 + \gamma + (n-i)\left[ \frac{(n-i-1)(1-\wtb)}{2} - 1 \right]. 
\end{align}

\subsection{Case 3: $0 < \wtb < \frac{n-3}{n-1}$} \label{App:2.2}
Note that, for $n = 3$, the interval is empty. So, we shall only consider $n \ge 4$.

Let $\wtb = \alpha \frac{n-3}{n-1}$, $\alpha \in (0,1)$. Then, $(n-1)(1-\wtb)-2 = (1-\alpha)(n-3) > 0$. Next, $(n-i)(1-\wtb)-1=0$ implies $i = n - \frac{1}{1-\wtb} = n-\gamma+1 =: \widetilde{i}$, because $\gamma = \frac{2-\wtb}{1-\wtb} = 1 + \frac{1}{1-\wtb}$. Thus, $(n-i)(1-\wtb) - 1 > 0$ when $i < \widetilde{i}$, and $(n-i)(1-\wtb) - 1 < 0$ when $i > \widetilde{i}$. Let $i^* = \lfloor \widetilde{i} \rfloor$. Then
\begin{align*}
    \begin{cases}
        (n-i)(1-\wtb)-1 \ge 0,&i \le i^*, \\
        (n-i)(1-\wtb)-1 < 0,&i > i^*.
    \end{cases}
\end{align*}
For $\wtb \in \left(0,\frac{n-3}{n-1}\right)$, we have the following facts: 
\begin{enumerate}
\item[(i)] $\gamma \in \left(2,\frac{n+1}{2}\right)$. Thus, $\gamma < \frac{n+1}{2}$; 
\item[(ii)] $i^* \in \left[ \frac{n}{2}, n-2 \right]$, because $\widetilde{i} = n-\gamma+1$, we have $\frac{n-1}{2} < n-\gamma < i^* \le n - \gamma + 1 < n-1$;
\item[(iii)] $n - \gamma \in \left(\frac{n-1}{2},n-2\right)$; 
\item[(iv)] $\lfloor \gamma \rfloor \in \left[2, \frac{n}{2} \right]$, therefore $\lfloor \gamma \rfloor \le i^*$. The equality only holds if $n$ is even.
\end{enumerate}

\paragraph{\underline{Subcase 3.1.1:} $1 \leq i < \lfloor \gamma \rfloor$.} Note that 
$$
0  \ge \lfloor \gamma \rfloor - \gamma > \dots > i - \gamma > \dots > 1 - \gamma
$$
and
$$
 (n-i)(1-\wtb)-1 > \dots > (n-i^*)(1-\wtb)-1 \ge 0 > (n-(i^*+1))(1-\wtb)-1 > -\wtb > -1.
$$
Thus, for $1 \leq i \leq \lfloor \gamma \rfloor$,
\begin{align}\label{eq:sumA31}
  \sum_{j=1}^i |j-\gamma| = \sum_{j=1}^i (\gamma - j) = \gamma i - \frac{1}{2} i(i+1)  
\end{align}
and
\begin{align}
    \sum_{j=i+1}^n |(n-j)(1-\wtb)-1| &= \sum_{j=i+1}^{i^*} |(n-j)(1-\wtb)-1| + \sum_{j=i^*+1}^n |(n-j)(1-\wtb)-1| \notag \\
    &= \sum_{j=i+1}^{i^*} ((n-j)(1-\wtb)-1) + \sum_{j=i^*+1}^n (1-(n-j)(1-\wtb)) \notag \\
    &\le  1 + \gamma + (n-i)\left[ \frac{(n-i-1)(1-\wtb)}{2} - 1 \right], \label{eq:sumA32}
\end{align}
from~\eqref{eq:App1.1} in Appendix~\ref{App:2.1}. With
\begin{align*}
    \left| \wtt_{i,1}^{-1} \right| = \left| (1-\gamma) \frac{(n-i)(1-\wtb)-1}{(n-1)(1-\wtb)-2} \right| = (\gamma-1) \frac{(n-i)(1-\wtb)-1}{(n-1)(1-\wtb)-2},
\end{align*}
and
\begin{align*}
    \left| \wtt_{i,n}^{-1} \right| = \left| \wtt_{n,i}^{-1} \right| = \left| \frac{i-\gamma}{(n-1)(1-\wtb)-2} \right| = \frac{\gamma - i}{(n-1)(1-\wtb)-2},
\end{align*}
the sum~\eqref{eq1_abs_sum} can be bounded from above as follows:
\begin{align*}
    \sum_{j=1}^n \left| \wtt_{i,j}^{-1} \right| &\le \frac{(n-i)(1-\wtb)-1}{(n-1)(1-\wtb)-2} \left( \gamma i - \frac{1}{2} i(i+1) \right) \\
    &+ \frac{\gamma - i}{(n-1)(1-\wtb)-2} \left( 1 + \gamma + (n-i)\left[ \frac{(n-i-1)(1-\wtb)}{2} - 1 \right] \right) \\
    &= - \left( \sum_{j=1}^n \wtt_{i,j}^{-1} + \frac{(i-\gamma)(1+\gamma)}{(n-1)(1-\wtb)-2} \right) = -P_i,
\end{align*}
where $P_i$ is defined in \eqref{eq:Pi}. In this case,
\begin{align*}
    P_i' = -\left(-1 + \frac{n+1}{2} + \frac{\gamma+1}{(n-1)(1-\wtb)-2} \right) < 0.
\end{align*}
$-P_i$ is then a decreasing function of $i$. It attains maximum at $i = 1$. This leads to the following bound:
\begin{align} \label{eq:app2.3.1}
  \max_{1 \le i < \lfloor \gamma \rfloor} \sum_{j=1}^n \left| \wtt_{i,j}^{-1} \right| \le -P = \max_{1 \le i < \lfloor \gamma \rfloor} -P_i = -P_1 = \frac{n(\gamma-1)}{2} + \frac{(\gamma-1)^2(\gamma+1)}{n+1-2\gamma}.
\end{align}

\paragraph{\underline{Subcase 3.1.2:} }
Note that when $\lfloor \gamma \rfloor < i^*$, the above analysis still applies with $i = \lfloor \gamma \rfloor$. In this scenario, we obtain the same upper bound, denoted as $-P$. However, this bound may not hold when $i = \lfloor \gamma \rfloor = i^* = \frac{n}{2}$ with $n$ is even. Therefore, we need to separately evaluate the $i^*$-th sum in this case and compare it with $-P_1$. We have
\begin{align*}
1-\gamma < 2-\gamma < \ldots <  i - \gamma = i^* - \gamma =  \lfloor \gamma \rfloor - \gamma  \leq 0
\end{align*}
and
\begin{align*}
0 > (n-i^*-1)(1-\wtb)-1 > \ldots > -1.
\end{align*}
Therefore,
\begin{align*}
\sum_{j=1}^{i} \left| j - \gamma \right| = \sum_{j=1}^{i}(\gamma - j) = \gamma i - \frac{i(i+1)}{2} = \frac{n\gamma}{2} - \frac{n(n+2)}{8},
\end{align*}
\begin{align*}
\sum_{j=i+1}^{n}\left| (n-j)(1-\wtb)-1\right| 
    &= \sum_{j=i+1}^{n}(1 - (n-j)(1-\wtb))\\ 
    &= (n-i) - \frac{(n-i)(n-i-1)(1-\wtb)}{2}\\
    &= \frac{n}{2} - \frac{n(n-2)(1-\wtb)}{8},
\end{align*}
and
\begin{align*}
\left|\wtt^{-1}_{i,1} \right| &= \frac{\gamma-1}{2} \frac{n(1-\wtb)-2}{(n-1)(1-\wtb)-2}\\
\left|\wtt^{-1}_{i,n} \right| &= \frac{\gamma - \frac{n}{2}}{(n-1)(1-\wtb)-2}
\end{align*}
Using this, we get the $i^*$-th sum
\begin{align*}
\sum_{j=1}^{n} \left|\wtt^{-1}_{i^*,j} \right| 
    &= \frac{1}{2}\frac{n(1-\wtb)-2}{(n-1)(1-\wtb)-2} \left[  \frac{n\gamma}{2} - \frac{n(n+2)}{8}\right] + \frac{\gamma - \frac{n}{2}}{(n-1)(1-\wtb)-2} \left[\frac{n}{2} - \frac{n(n-2)(1-\wtb)}{8} \right]\\
    &= \frac{1}{16}\frac{n}{(n-1)(1-\wtb)-2}\left[(n(1-\wtb)-2)(4\gamma - n-2)+(2\gamma - n)(4 - (n-2)(1-\wtb))\right]\\
    &= \frac{1}{8} \frac{n}{(n-1)(1-\wtb)-2} \left[(1-\wtb)(\gamma n -2n +2\gamma) - n +2 \right]\\
    &= \frac{n}{8} \frac{\gamma n - 2n + 2\gamma}{n+1-2\gamma} - \frac{n(n-2)}{8} \frac{\gamma-1}{n+1-2\gamma}\\
    &= \frac{n}{8}  \frac{4\gamma - n - 2}{n+1-2\gamma}\\
    &= \frac{n(2\gamma-1)}{8(n+1-2\gamma)}- \frac{n}{8}.
\end{align*}
Now, let us compare this sum to $-P$ for $n\geq 4$
\begin{align*}
(-P) - \sum_{j=1}^{n} \left|\wtt^{-1}_{i^*,j} \right| 
    &= \frac{n(\gamma-1)}{2} + \frac{(\gamma-1)^2 (\gamma+1)}{n+1-2\gamma} -  \frac{n(2\gamma-1)}{8(n+1-2\gamma)} + \frac{n}{8}\\
    &= \frac{n(\gamma-1)}{2}+\frac{8(\gamma-1)^2 (\gamma+1) - n(4\gamma-n-2)}{8(n+1-2\gamma)} \\
    &\geq \frac{n(\gamma-1)}{2} + \frac{(n-2)^2 (n+2) - n^2}{8(n+1-2\gamma)}\\
    &\geq \frac{n(\gamma-1)}{2} + \frac{n^3 - 3n^2 -4n+8}{8(n+1-2\gamma)} \\
    &\geq \frac{n(\gamma-1)}{2} + \frac{(n-2)^2+4}{8(n+1-2\gamma)}>0.
\end{align*}
Therefore, we can say
\begin{align}
\max_{1\leq i \leq \lfloor \gamma \rfloor} \sum_{j=1}^{n} \left|\wtt^{-1}_{i,j} \right| \leq -P = \frac{n(\gamma-1)}{2} + \frac{(\gamma-1)^2 (\gamma+1)}{n+1-2\gamma}.
\end{align}

\paragraph{\underline{Subcase 3.2.1:} $\lfloor \gamma \rfloor < i < \frac{n+1}{2}.$}

By using the nested inequalities above~\eqref{eq:sumA31},
\begin{align}
    \sum_{j=1}^i |j-\gamma| &= \sum_{j=1}^{\lfloor \gamma \rfloor} |j-\gamma| + \sum_{j=\lfloor \gamma \rfloor + 1}^{i} |j-\gamma| = \sum_{j=1}^{\lfloor \gamma \rfloor} (\gamma - j) + \sum_{j=\lfloor \gamma \rfloor + 1}^{i} (j-\gamma) \notag \\
    &=\gamma \lfloor \gamma \rfloor - \frac{1}{2}\lfloor \gamma \rfloor (\lfloor \gamma \rfloor+1) + \frac{1}{2}i(i+1) - \gamma i - \frac{1}{2}\lfloor \gamma \rfloor (\lfloor \gamma \rfloor+1) + \gamma \lfloor \gamma \rfloor \notag \\
    &= \frac{1}{2}i(i+1) - \gamma i + \gamma^2 - \lfloor \gamma \rfloor  - (\gamma - \lfloor \gamma \rfloor)^2 \notag \\
    &\le \frac{1}{2}i(i+1) - \gamma i + \gamma^2 - \lfloor \gamma \rfloor  \notag \\
    &< \frac{1}{2}i(i+1) - \gamma i + \gamma^2 -\gamma + 1. \label{eq:sumA33}
\end{align}
where the first and the second inequality are obtained after applying $0 \le \gamma - \lfloor \gamma \rfloor < 1$. Also by \eqref{eq:App1.1},
\begin{align}
    \sum_{j=i+1}^n |(n-j)(1-\wtb)-1| &= \sum_{j=i+1}^{i^*} |(n-j)(1-\wtb)-1| + \sum_{j=i^*+1}^n |(n-j)(1-\wtb)-1| \notag \\
    &= \sum_{j=i+1}^{i^*} ((n-j)(1-\wtb)-1) + \sum_{j=i^*+1}^n (1-(n-j)(1-\wtb)) \notag \\
    &\le 1 + \gamma + (n-i)\left[ \frac{(n-i-1)(1-\wtb)}{2} - 1 \right]. \label{eq:sumA34}
\end{align}
With 
\begin{align*}
    \left| \wtt_{i,1}^{-1} \right| = \left| (1-\gamma) \frac{(n-i)(1-\wtb)-1}{(n-1)(1-\wtb)-2} \right| = (\gamma-1) \frac{(n-i)(1-\wtb)-1}{(n-1)(1-\wtb)-2},
\end{align*}
and
\begin{align*}
    \left| \wtt_{i,n}^{-1} \right| = \left| \wtt_{n,i}^{-1} \right| = \left| \frac{i-\gamma}{(n-1)(1-\wtb)-2} \right| = \frac{i - \gamma}{(n-1)(1-\wtb)-2},
\end{align*}
the sum~\eqref{eq1_abs_sum} can be bounded from above as follows, after simplification,
\begin{align*}
    \sum_{j=1}^n \left| \wtt_{i,j}^{-1} \right| &< \underbrace{\sum_{j=1}^n \wtt^{-1}_{i,j} + \frac{(i-\gamma)(1+\gamma)}{(n-1)(1-\wtb)-2}}_{P_i} + \frac{(n-i)(1-\wtb)-1}{(n-1)(1-\wtb)-2} (\gamma^2 - \gamma + 1) =: M_i.
\end{align*}
For $\lfloor \gamma \rfloor < i <  \frac{n+1}{2}$,
\begin{align*}
    M_i' = -i + \frac{n+1}{2} + \frac{\gamma+1}{(n-1)(1-\wtb)-2} - \frac{(1-\wtb)(\gamma^2-\gamma+1)}{(n-1)(1-\wtb)-2} = -i + \frac{n+1}{2} + \frac{\gamma-2}{n+1-2\gamma} > 0,
\end{align*}
because of Fact (i): $\gamma \in \left(2,\frac{n+1}{2}\right)$. So, $M$ is an increasing function in the interval $\left[1,  \frac{n+1}{2} \right]$, with the maximum attained at $i = \frac{n+1}{2} + \frac{\gamma - 2}{n+1-2\gamma} > \frac{n+1}{2}$.
Therefore,
\begin{align} \label{eq:app2.3.2}
    \max_{i \in  (\lfloor \gamma \rfloor, \frac{n+1}{2}}) \sum_{j=1}^n \left| \wtt_{i,j}^{-1} \right| &< M_{\lfloor \frac{n+1}{2} \rfloor} \le R: = M_{\frac{n+1}{2}}  = \frac{(n+1)^2}{8} - \gamma \left( \frac{n+1}{2} - \gamma \right).
\end{align}

\paragraph{\underline{Subcase 3.2.2:}}
Note that for \( i^* > \frac{n+1}{2} \), the above analysis still holds with \( i = \frac{n+1}{2} \). In that case, we obtain the same upper bound, which is \( R \). However, this bound might not hold when \( i = i^* = \frac{n+1}{2} \) with \( n \) being odd. Therefore, we need to consider this case separately. Note \( \frac{n-1}{2} < \gamma < \frac{n+1}{2} \), which implies \( \lfloor \gamma \rfloor = \frac{n-1}{2} \). Since we already know the upper bound for \( \lfloor \gamma \rfloor < i < \frac{n+1}{2} = i^* \), we need to find the sum for \( i = i^* = \frac{n+1}{2} \).
We can observe the following inequalities
\begin{align*}
    1-\gamma < \ldots < \lfloor \gamma \rfloor - \gamma \leq 0 < i - \gamma = i^* - \gamma = \frac{n+1}{2} - \gamma 
\end{align*}
and
\begin{align*}
    0 > (n-i^*-1)(1-\wtb)-1 > \ldots > -1.
\end{align*}
The corresponding sums are calculated as follows.
For the first sum
\begin{align*}
    \sum_{j=1}^{i}\left|j-\gamma \right| 
        &= \sum_{j=1}^{\lfloor \gamma \rfloor} (\gamma - j) + \frac{n+1}{2}-\gamma\\
        &= \gamma \lfloor \gamma \rfloor - \frac{\lfloor\gamma \rfloor^2 + \lfloor \gamma \rfloor}{2} + \frac{n+1}{2} -\gamma\\
        &= \gamma  \frac{n-1}{2} - \frac{(n-1)^2}{8} - \frac{n-1}{4} + \frac{n+1}{2} - \gamma\\
        &= \gamma  \frac{n-3}{2} -\frac{(n-1)^2}{8} + \frac{n+3}{4}.
\end{align*}
For the second sum
\begin{align*}
    \sum_{j=i+1}^{n} \left|(n-j)(1-\wtb)-1\right| 
        &= (n-i)\left[1 - \frac{(n-i-1)(1-\wtb)}{2}\right] =\frac{n-1}{2} \left[1 - \frac{(n-3)(1-\wtb)}{4} \right].
\end{align*}
We also have
\begin{align*}
    \left| \wtt^{-1}_{i,1} \right| = \left| \wtt^{-1}_{i,n} \right|=  \frac{\gamma - 1}{2}.
\end{align*}
Substituting these values into the formula, we get
\begin{align*}
    \sum_{j=1}^{n} \left| \wtt^{-1}_{i^*,j} \right| 
        &= \frac{1}{2} \left( \gamma  \frac{n-3}{2} - \frac{(n-1)^2}{8} + \frac{n+3}{4} \right) + \frac{(\gamma-1)(n-1)}{4} \left( 1 - \frac{(n-3)(1-\wtb)}{4}\right)\\
        &=\frac{\gamma(n-2)}{2} - \frac{(n-3)(n+1)}{8}.
\end{align*}
Now, let us compare this sum to the above upper bound $R$ for $n\geq 4$
\begin{align*}
    \frac{(n+1)^2}{8} - \gamma\left(\frac{n+1}{2} - \gamma \right) - \frac{\gamma(n-2)}{2} + \frac{(n-3)(n+1)}{8} &= \frac{n^2-1}{4} - \frac{\gamma(2n-1)}{2} + \gamma^2\\
    &= \left(\gamma - \frac{2n-1}{4}\right)^2 + \frac{4n-5}{16} >0.
\end{align*}
Therefore, we can state the following upper bound
\begin{align*}
    \max_{\lfloor \gamma \rfloor < i \leq \frac{n+1}{2}} \sum_{j=1}^{n}\left| \wtt^{-1}_{i,j}\right| \leq \frac{(n+1)^2}{8} - \gamma \left(\frac{n+1}{2} - \gamma \right).
\end{align*}
Combining the results from Subsections 3.1.2, 3.2.2 and using the centrosymmetricity of $\wtT^{-1}_{n}$, we get 
\begin{align}
    \|\wtT^{-1}_{n}\|_{\infty} = \max_{1 \le i \le \frac{n+1}{2}} \sum_{j=1}^{n} |\wtt^{-1}_{i,j}|\le \max \left\{\frac{n(\gamma-1)}{2}+\frac{(\gamma-1)^2 (\gamma+1)}{n+1-2\gamma}, \frac{(n+1)^2}{8} - \gamma\left(\frac{n+1}{2} - \gamma \right)\right\}.
\end{align}

\pagebreak

\end{document}